\documentclass[10pt]{amsart}
\usepackage{amsfonts}
\usepackage{amssymb}
\usepackage{amsfonts}
\usepackage{bbm}
\newtheorem{thm1}{Theorem}[section]
\newtheorem{lem1}[thm1]{Lemma}
\newtheorem{rem1}[thm1]{Remark}
\newtheorem{def1}[thm1]{Definition}

\newtheorem{cor1}[thm1]{Corollary}

\newtheorem{prop1}[thm1]{Proposition}
\newtheorem{ex1}[thm1]{Example}
\newtheorem{not1}[thm1]{Notation}
\begin{document}

\title[]
{On the binomial arithmetical rank of lattice ideals}
\author[A. Katsabekis]{Anargyros Katsabekis}
\address {Centrum Wiskunde \& Informatica (CWI), Postbus 94079, 1090 GB Amsterdam, The Netherlands} \email{katsabek@aegean.gr}
\thanks{This work was carried out during the tenure of an ERCIM \textit{Alain Bensoussan} Fellowship Programme. The research leading to these results has received funding from the European Union Seventh Framework Programme (FP7/2007-2013) under grant agreement no 246016.}


\begin{abstract} To any lattice $L \subset \mathbb{Z}^{m}$ one can associate the lattice ideal $I_{L} \subset K[x_{1},\ldots,x_{m}]$. This paper concerns the study of the relation between the binomial arithmetical rank and the minimal number of generators of $I_{L}$. We provide lower bounds for the binomial arithmetical rank and the $\mathcal{A}$-homogeneous arithmetical rank of $I_{L}$. Furthermore, in certain cases we show that the binomial arithmetical rank equals the minimal number of generators of $I_{L}$. Finally we consider a class of determinantal lattice ideals and study some algebraic properties of them.
\end{abstract}
\maketitle


\section{Introduction}
\par Let $K[x_{1},\ldots,x_{m}]$ be a polynomial ring in $m$ variables over any field $K$. As usual, we will denote by ${\bf x}^{\bf u}$ the monomial $x_{1}^{u_{1}} \cdots x_{m}^{u_m}$ of $K[x_{1},\ldots,x_{m}]$, with ${\bf u}=(u_{1},\ldots,u_{m}) \in \mathbb{N}^{m}$, where $\mathbb{N}$ stands for the set of non-negative integers. A {\em binomial} is a polynomial which is a difference of two monomials. A {\em binomial ideal} is an ideal generated by binomials. Recall that a lattice is a finitely generated free abelian group. Given a lattice $L \subset \mathbb{Z}^{m}$, the ideal $$I_{L}=(\{{\bf x}^{{\bf u}_+}-{\bf x}^{{\bf u}_-}|{\bf u}={\bf u}_{+}-{\bf u}_{-} \in L\}) \subset K[x_{1},\ldots,x_{m}]$$ is called {\em lattice ideal}. Here ${\bf u}_{+} \in \mathbb{N}^{m}$ and ${\bf u}_{-} \in \mathbb{N}^{m}$ denote the positive and negative part of ${\bf u}$, respectively.

Throughout this paper we assume that $L$ is a non-zero positive sublattice of $\mathbb{Z}^{m}$, that is
$L \cap \mathbb{N}^{m}=\{\bf 0\}$. By the graded Nakayama's Lemma, all minimal binomial generating sets of $I_{L}$ have the same cardinality. The cardinality of any minimal generating set of $I_{L}$ consisting of binomials is commonly known as the {\em minimal number of generators} of $I_{L}$, denoted by $\mu(I_{L})$.

If $L=<{\bf l}_{1},\ldots,{\bf l}_{k}>$ is a sublattice of $\mathbb{Z}^{m}$ of rank $k<m$, then the {\em saturation} of $L$ is the lattice $$Sat(L):=\{{\bf u} \in \mathbb{Z}^{m}|d{\bf u} \in L \ \textrm{for some} \ d \in \mathbb{Z}, \ d \neq 0\}.$$ Clearly, the inclusion $L \subset Sat(L)$ holds. Also there exists a set of vectors $\mathcal{A}=\{{\bf a}_{1},\ldots,{\bf a}_{m}\} \subset \mathbb{Z}^{n}$ such that $Sat(L)=ker_{\mathbb{Z}}(\mathcal{A})$, where $n=m-k$ and $$ker_{\mathbb{Z}}(\mathcal{A}):=\{(q_{1},\ldots,q_{m}) \in \mathbb{Z}^{m}|q_{1}{\bf a}_{1}+\cdots+q_{m}{\bf a}_{m}={\bf 0}\}.$$
When $L$ is saturated, i.e. $L=Sat(L)$, the ideal $I_{L}$ is called {\em toric ideal}. We will write for simplicity $I_{\mathcal{A}}:=I_{ker_{\mathbb{Z}}(\mathcal{A})}$. The toric ideal $I_{\mathcal{A}}$ is the kernel of the $K$-algebra homomorphism $\phi: K[x_{1},\ldots,x_{m}] \rightarrow K[t_{1},\ldots,t_{n}]$ given by $\phi(x_{i})={\bf t}^{{\bf a}_{i}}$, for every $i=1,\ldots,m$ (see \cite{St}). Thus every toric ideal is prime.

 We grade $K[x_{1},\ldots,x_{m}]$ by setting ${\rm deg}_{\mathcal{A}}(x_{i})={\bf a}_{i}$, $1 \leq i \leq m$. The $\mathcal{A}$-degree of the monomial ${\bf x}^{\bf u}$ is ${\rm deg}_{\mathcal{A}}({\bf x}^{\bf u})=u_{1}{\bf a}_{1}+\cdots+u_{m}{\bf a}_{m} \in \mathbb{N}\mathcal{A}$ where $\mathbb{N}\mathcal{A}$ is the semigroup generated by $\mathcal{A}$. A polynomial $F$ is called $\mathcal{A}$-{\em homogeneous} if the monomials in each nonzero term of $F$ have the same $\mathcal{A}$-degree. The ideal $I$ is called $\mathcal{A}$-homogeneous if it is generated by $\mathcal{A}$-homogeneous polynomials. The lattice ideal $I_{L}$ is $\mathcal{A}$-homogeneous, since it is generated by binomials and every binomial ${\bf x}^{{\bf u}_{+}}-{\bf x}^{{\bf u}_{-}}$ is $\mathcal{A}$-homogeneous.

For an ideal $I \subset K[x_{1},\ldots,x_{m}]$ we let $rad(I)$ be its radical. The {\em arithmetical rank} of $I_{L}$, denoted by ${\rm ara}(I_{L})$, is the smallest integer $s$ for which there exist polynomials $F_{1},\ldots,F_{s}$ in $I_L$ such that $rad(I_{L})=rad(F_{1},\ldots,F_{s})$. When $K$ is algebraically closed, ${\rm ara}(I_{L})$ is the smallest number of hypersurfaces whose intersection is set-theoretically equal to the algebraic set defined by $I_{L}$. Computing the arithmetical rank is one of the classical problems of Algebraic Geometry which remains open even for very simple cases, like the ideal of the Macaulay curve in the three-dimensional projective space. If all the polynomials $F_{1},\ldots,F_{s}$ satisfying $rad(I_{L})=rad(F_{1},\ldots,F_{s})$ are $\mathcal{A}$-homogeneous, the smallest integer $s$ is called the $\mathcal{A}$-{\em homogeneous arithmetical rank} of $I_{L}$ and will be denoted by ${\rm ara}_{\mathcal{A}}({I_{L}})$. Since $I_L$ is generated by binomials, it is natural to define the {\em binomial arithmetical rank} of $I_{L}$, denoted by ${\rm bar}(I_L)$, as the smallest integer $s$ for which there exist binomials $B_{1},\ldots,B_{s}$ in $I_L$ such that $rad(I_{L})=rad(B_{1},\ldots,B_{s})$. From the definitions and the generalized Krull's principal ideal theorem we have the following inequalities: $${\rm ht}(I_{L}) \leq {\rm ara}(I_{L}) \leq {\rm ara}_{\mathcal{A}}(I_{L}) \leq {\rm bar}(I_{L}) \leq \mu(I_{L}).$$
Where ${\rm ht}(I_{L})$ is the height of $I_{L}$ which equals the rank of the lattice $L$, see Corollary 2.2 on \cite{E-S}.

In this paper we are interested in the problem when the equality ${\rm bar}(I_{L})=\mu(I_{L})$ holds. Clearly it is valid for the special class of complete intersection lattice ideals. Recall that a lattice ideal $I_{L}$ is {\em complete intersection} if $\mu(I_{L})={\rm ht}(I_{L})$. The above problem was considered for the case of toric ideals associated with finite graphs in \cite{K}, see section 3 for the definition of such ideals. More precisely the author reveals two cases in which the binomial arithmetical rank coincides with the minimal number of generators for the toric ideal $I_{\mathcal{A}_{G}}$ of a graph $G$, namely when $G$ is bipartite or the ideal $I_{\mathcal{A}_G}$ is generated by quadratic binomials. The main aim of this work is to generate new classes of lattice ideals for which the equality ${\rm bar}(I_{L})=\mu(I_{L})$ holds.

In section 2 we consider the indispensable monomials of a lattice ideal $I_{L}$ and study the related simplicial complex $\Gamma_{L}$. We provide a necessary condition for the generation of the radical of a lattice ideal $I_{L}$ up to radical, see Theorem 2.9. Using this result and also the notion of $J$-matchings in simplicial complexes, introduced in \cite{KT}, we obtain lower bounds for the binomial arithmetical rank and the $\mathcal{A}$-homogeneous arithmetical rank of a lattice ideal (see Theorem 2.13), which are in general different (see the discussion after Theorem 2.13 and also Example 2.14) than the bounds given in Theorem 5.6 of \cite{KMT}.

In section 3 we deal with the toric ideals associated with graphs. After presenting the basic theory of such ideals, we concentrate ourselves on the case that the graph satisfies a certain condition, which guarantees that the toric ideal is generated by binomials of a specific form. We use Theorem 2.13 to show that the equality ${\rm bar}(I_{\mathcal{A}_{G}})=\mu(I_{\mathcal{A}_{G}})$ holds under a mild assumption on the toric ideal $I_{\mathcal{A}_{G}}$ (see Theorem 3.14). This assumption is fulfilled by the toric ideal associated with a bipartite graph, as well as a toric ideal generated by quadratic binomials. As applications we prove that the binomial arithmetical rank equals the minimal number of generators of $I_{\mathcal{A}_{G}}$ for two types of graphs, namely the wheel graph and a weakly chordal graph.

Section 4 is devoted to the study of a class of determinantal ideals $I_{2}(D)$ with the property that ${\rm bar}(I_{2}(D))=\mu(I_{2}(D))$. Every such ideal is a lattice ideal, so it is of the form $I_{L}$ for a certain lattice $L$, and also has a unique minimal system of binomial generators. Finally we consider the lattice basis ideal $J_{L}$ and determine its minimal primary decomposition, under the condition that the ideal $I_{2}(D)$ is prime.



\section{General results}

Let $L \subset \mathbb{Z}^{m}$ be a non-zero positive lattice with $Sat(L)=ker_{\mathbb{Z}}(\mathcal{A})$, where $\mathcal{A}=\{{\bf a}_{1},\ldots,{\bf a}_{m}\} \subset \mathbb{Z}^{n}$. In this section we associate to $L$ the simplicial complex $\Gamma_{L}$. We show that combinatorial invariants of this complex provide lower bounds for the binomial arithmetical rank and the $\mathcal{A}$-homogeneous arithmetical rank of $I_{L}$.

\begin{not1} {\rm For a vector ${\bf v}=(v_{1},\ldots,v_{m}) \in \mathbb{Z}^{m}$, we shall denote by ${\rm supp}({\bf v}):=\{i \in \{1,\ldots,m\}|v_{i} \neq 0\}$ the support of ${\bf v}$. Given a monomial ${\bf x}^{\bf w} \in K[x_{1},\ldots,x_{m}]$, we let ${\rm supp}({\bf x}^{\bf w})={\rm supp}({\bf w})$.}
\end{not1}

\begin{def1} \label{Indisp} {\rm A {\em binomial} $B=M-N \in I_{L}$ is called {\em indispensable} of $I_{L}$ if every system of binomial generators of $I_{L}$ contains $B$ or $-B,$ while a {\em monomial} $M$ is called {\em indispensable} of $I_{L}$ if every system of binomial generators of $I_{L}$ contains a binomial $B$ such that $M$ is a monomial of $B$.}
\end{def1}
Let $\mathcal{M}_{L}$ be the ideal generated by all monomials $M$ for which there exists a nonzero $M-N \in I_{L}$. Proposition 1.5 of \cite{KO} implies that the set of indispensable monomials of $I_{L}$ is the unique minimal generating set of $\mathcal{M}_{L}$.

\begin{rem1} \label{RemIndi} {\rm If $\{B_{1}=M_{1}-N_{1},\ldots,B_{s}=M_{s}-N_{s}\}$ is a generating set of $I_{L}$, then $\mathcal{M}_{L}=(M_{1},\ldots,M_{s},N_{1},\ldots,N_{s})$.}
\end{rem1}

Let $\mathcal{T}$ be the set of all $E \subset \{1,\ldots,m\}$ such that $E={\rm supp}(M)$, where $M$ is an indispensable monomial of $I_{L}$. We shall denote by $\mathcal{T}_{\rm min}$ the set of minimal elements of $\mathcal{T}$.

\begin{def1} {\rm  We associate to $L$ the simplicial complex $\Gamma_{L}$ with vertices the elements of $\mathcal{T}_{\rm min}$. Let $T=\{E_{1},\ldots,E_{k}\}$ be a subset of $\mathcal{T}_{\rm min}$, then $T \in \Gamma_{L}$ if \begin{enumerate} \item for every $E_{i}$, $1 \leq i \leq k$, there exists a monomial $M_{i}$ with ${\rm supp}(M_{i})=E_{i}$ and \item the monomials $M_{1},\ldots,M_{k}$ have the same $\mathcal{A}$-degree, i.e. it holds that $${\rm deg}_{\mathcal{A}}(M_{1})={\rm deg}_{\mathcal{A}}(M_{2})=\cdots={\rm deg}_{\mathcal{A}}(M_{k}).$$
\end{enumerate}}
\end{def1}

A non-zero vector ${\bf u}=(u_{1},\ldots,u_{m}) \in ker_{\mathbb{Z}}(\mathcal{A})$ is called a {\em circuit} if its support is minimal with respect to inclusion, namely there exists no other vector ${\bf v} \in ker_{\mathbb{Z}}(\mathcal{A})$ such that ${\rm supp}({\bf v}) \subsetneqq {\rm supp}({\bf u})$, and the coordinates of ${\bf u}$ are relatively prime. The binomial ${\bf x}^{{\bf u}_{+}}-{\bf x}^{{\bf u}_{-}} \in I_{\mathcal{A}}$ is called also circuit. We will make the connection between the elements of $\Gamma_{L}$ and the circuits of $I_{\mathcal{A}}$.

\begin{lem1} \label{Nosubset} If $E \in \mathcal{T}_{\rm min}$, then \begin{enumerate} \item there exists no circuit ${\bf x}^{{\bf u}_{+}}-{\bf x}^{{\bf u}_{-}} \in I_{\mathcal{A}}$ such that ${\rm supp}({\bf x}^{{\bf u}_{+}}) \subsetneqq E$ or ${\rm supp}({\bf x}^{{\bf u}_{-}}) \subsetneqq E$. \item there exists a circuit ${\bf x}^{{\bf u}_{+}}-{\bf x}^{{\bf u}_{-}} \in I_{\mathcal{A}}$ such that ${\rm supp}({\bf x}^{{\bf u}_{+}})=E$ or ${\rm supp}({\bf x}^{{\bf u}_{-}})=E$.
\end{enumerate}
\end{lem1}
\noindent \textbf{Proof.} (1) Suppose that $I_{\mathcal{A}}$ has a circuit ${\bf x}^{{\bf u}_{+}}-{\bf x}^{{\bf u}_{-}}$ such that ${\rm supp}({\bf x}^{{\bf u}_{+}}) \subsetneqq E$. Since ${\bf u} \in ker_{\mathbb{Z}}(\mathcal{A})=Sat(L)$, there exists a positive integer $d$ such that $L \ni d{\bf u}={\bf v}$. Notice that ${\rm supp}({\bf v}_{+})={\rm supp}({\bf u}_{+})$ and ${\rm supp}({\bf v}_{-})={\rm supp}({\bf u}_{-})$. Since ${\rm supp}({\bf v}_{+}) \subsetneqq E$ and $E \in \mathcal{T}_{\rm min}$, the monomial ${\bf x}^{{\bf v}_{+}}$ is not indispensable. Thus there exists an indispensable monomial $M$ of $I_{L}$ such that $M$ divides ${\bf x}^{{\bf v}_{+}}$ and $M \neq {\bf x}^{{\bf v}_{+}}$. As a consequence ${\rm supp}(M) \subseteq {\rm supp}({\bf x}^{{\bf v}_{+}}))$ and therefore ${\rm supp}(M) \subsetneqq E$, a contradiction to the fact that $E \in \mathcal{T}_{\rm min}$.\\ (2) Let $E={\rm supp}({\bf x}^{{\bf v}_{+}})$ where ${\bf x}^{{\bf v}_{+}}-{\bf x}^{{\bf v}_{-}} \in I_{L}$ and ${\bf x}^{{\bf v}_{+}}$ is an indispensable monomial of $I_{L}$. The vector ${\bf v}={\bf v}_{+}-{\bf v}_{-} \in ker_{\mathbb{Z}}(\mathcal{A})$ and therefore there exists, from Proposition 4.10 of \cite{St}, a circuit ${\bf u}$ conformal to ${\bf v}$, i.e. ${\rm supp}({\bf u}_{+}) \subseteq {\rm supp}({\bf v}_{+})$ and ${\rm supp}({\bf u}_{-}) \subseteq {\rm supp}({\bf v}_{-})$. Thus ${\rm supp}({\bf x}^{{\bf u}_{+}}) \subseteq {\rm supp}({\bf x}^{{\bf v}_{+}})=E$, so we have, from (1), that necessarily $E={\rm supp}({\bf x}^{{\bf u}_{+}})$. \hfill $\square$\\

We shall denote by $\mathcal{C}_{\mathcal{A}}$ the set of circuits of $\mathcal{A}$. Put $$\mathcal{C}:=\{E \subset \{1,\ldots,m\}| {\rm supp}({\bf u}_{+})=E \ \textrm{or} \ {\rm supp}({\bf u}_{-})=E \ \textrm{where} \ {\bf u} \in \mathcal{C}_{\mathcal{A}}\}.$$ Let $\mathcal{C}_{\rm min}$ be the set of minimal elements of $\mathcal{C}$.

\begin{prop1} \label{Cmin} It holds that $\mathcal{T}_{\rm min}=\mathcal{C}_{\rm min}$.
\end{prop1}
\noindent \textbf{Proof.} By Lemma \ref{Nosubset} we have that $\mathcal{T}_{\rm min} \subseteq \mathcal{C}_{\rm min}$. Conversely consider a set $E \in \mathcal{C}_{\rm min}$, then $E={\rm supp}({\bf x}^{{\bf u}_{+}})$ where ${\bf x}^{{\bf u}_{+}}- \in {\bf x}^{{\bf u}_{-}} \in I_{\mathcal{A}}$ is a circuit. Since ${\bf u} \in ker_{\mathbb{Z}}(\mathcal{A})=Sat(L)$, there exists a positive integer $d$ such that $L \ni d{\bf u}={\bf v}$. Notice that ${\rm supp}({\bf x}^{{\bf v}_{+}})={\rm supp}({\bf x}^{{\bf u}_{+}})$ and ${\rm supp}({\bf x}^{{\bf v}_{-}})={\rm supp}({\bf x}^{{\bf u}_{-}})$. Since ${\bf x}^{{\bf v}_{+}}$ belongs to the monomial ideal $\mathcal{M}_{L}$, there exists an indispensable monomial $M$ of $I_{L}$ with ${\rm supp}(M) \in \mathcal{T}_{\rm min}$ such that ${\rm supp}(M) \subseteq {\rm supp}({\bf x}^{{\bf v}_{+}})$. Now Lemma \ref{Nosubset} implies that ${\rm supp}(M) \in \mathcal{C}_{\rm min}$. But ${\rm supp}(M) \subseteq E$ and also $E \in \mathcal{C}_{\rm min}$, so $E={\rm supp}(M)$. \hfill $\square$\\

\begin{rem1} \label{RemCircuE} {\rm (1) In \cite{KT} a simplicial complex $\Delta_{\mathcal{A}}$ is associated to the vector configuration $\mathcal{A}$. By Proposition \ref{Cmin} the simplicial complex $\Gamma_{L}$ has the same vertex set with $\Delta_{\mathcal{A}}$. It is not hard to check that they are actually identical.\\ (2) By Theorem 4.2 (ii) of \cite{KT}, $\{E,E'\}$ is an edge of $\Gamma_{L}$ if and only if there is a circuit ${\bf x}^{{\bf u}_{+}}-{\bf x}^{{\bf u}_{-}} \in I_{\mathcal{A}}$ such that $E={\rm supp}({\bf x}^{{\bf u}_{+}})$ and $E'={\rm supp}({\bf x}^{{\bf u}_{-}})$.}
\end{rem1}

\begin{ex1} \label{BasicExample} {\rm Consider the lattice $L=ker_{\mathbb{Z}}(\mathcal{A})$ where $\mathcal{A}$ is the set of columns of the matrix $$P=\begin{pmatrix}
1 & 1 & 1 & 0 & 0 & 0 & 0 & 0 & 0 & 0 & 0 & 0
            \\
            0 & 0 & 0 & 1 & 1 & 1 & 3 & 0 & 1 & 2 & 1 & 2
            \\
            1 & 1 & 0 & 1 & 0 & 0 & 0 & 0 & 1 & 0 & 1 & 2
            \\
            1 & 0 & 1 & 0 & 1 & 0 & 1 & 0 & 1 & 1 & 1 & 2
            \\
            0 & 1 & 1 & 1 & 1 & 2 & 1 & 2 & 0 & 1 & 1 & 0
            \\
            0 & 0 & 0 & 0 & 0 & 0 & 0 & 2 & 0 & 0 & 1 & 0
            \\
            0 & 0 & 0 & 0 & 0 & 0 & 0 & 2 & 0 & 0 & 1 & 0

            \end{pmatrix}.$$ The toric ideal $I_{\mathcal{A}}$ is minimally generated by the following binomials:\\ $B_{1}=x_{2}x_{5}-x_{3}x_{4}$, $B_{2}=x_{1}x_{6}-x_{3}x_{4}$, $B_{3}=x_{1}x_{4}-x_{2}x_{9}$, $B_{4}=x_{1}x_{5}-x_{3}x_{9}$, $B_{5}=x_{4}x_{5}-x_{6}x_{9}$, $B_{6}=x_{10}^2-x_{5}x_{7}$, $B_{7}=x_{11}^{2}-x_{8}x_{9}^{2}$, $B_{8}=x_{9}^{2}-x_{12}$.\\ The circuits of $\mathcal{A}$ are $$\mathcal{C}_{\mathcal{A}}=\{x_{2}x_{5}-x_{3}x_{4},x_{1}x_{6}-x_{3}x_{4},x_{1}x_{6}-x_{2}x_{5},x_{1}x_{4}-x_{2}x_{9},x_{1}^{2}x_{4}^{2}-x_{2}^{2}x_{12},x_{1}x_{5}-x_{3}x_{9},$$ $$x_{1}^{2}x_{5}^{2}-x_{3}^{2}x_{12}, x_{4}x_{5}-x_{6}x_{9},x_{4}^{2}x_{5}^{2}-x_{6}^{2}x_{12},x_{10}^2-x_{5}x_{7},x_{11}^{2}-x_{8}x_{9}^{2},x_{11}^{2}-x_{8}x_{12},x_{9}^{2}-x_{12},$$ $$x_{2}x_{10}^2-x_{3}x_{4}x_{7},x_{2}x_{10}^2-x_{1}x_{6}x_{7}, x_{1}^{2}x_{6}-x_{2}x_{3}x_{9},x_{1}^{4}x_{6}^{2}-x_{2}^{2}x_{3}^{2}x_{12},x_{3}x_{4}^{2}-x_{2}x_{6}x_{9},$$ $$x_{3}^{2}x_{4}^{4}-x_{2}^{2}x_{6}^{2}x_{12}, x_{2}x_{5}^{2}-x_{3}x_{6}x_{9},x_{2}^{2}x_{5}^{4}-x_{3}^{2}x_{6}^{2}x_{12}, x_{1}x_{10}^{2}-x_{3}x_{7}x_{9},x_{1}^{2}x_{10}^{4}-x_{3}^{2}x_{7}^{2}x_{12},$$ $$x_{4}x_{10}^{2}-x_{6}x_{7}x_{9},x_{4}^{2}x_{10}^{4}-x_{6}^{2}x_{7}^{2}x_{12}, x_{2}^{2}x_{11}^{2}-x_{1}^{2}x_{4}^{2}x_{8},x_{3}^{2}x_{11}^{2}-x_{1}^{2}x_{5}^{2}x_{8},$$ $$x_{6}^{2}x_{11}^{2}-x_{4}^{2}x_{5}^{2}x_{8}, x_{2}x_{10}^{4}-x_{3}x_{6}x_{7}^{2}x_{9},x_{2}^{2}x_{10}^{8}-x_{3}^{2}x_{6}^{2}x_{7}^{4}x_{12}, x_{1}^{4}x_{6}^{2}x_{8}-x_{2}^{2}x_{3}^{2}x_{11}^{2},$$ $$x_{2}^{2}x_{6}^{2}x_{11}^{2}-x_{3}^{2}x_{4}^{4}x_{8}, x_{3}^{2}x_{6}^{2}x_{11}^{2}-x_{2}^{2}x_{5}^{4}x_{8},x_{3}^{2}x_{7}^{2}x_{11}^{2}-x_{1}^{2}x_{8}x_{10}^{4},x_{6}^{2}x_{7}^{2}x_{11}^{2}-x_{4}^{2}x_{8}x_{10}^{4}\}.$$ Thus the complex $\Gamma_{L}$ has 11 vertices: $E_{1}=\{1,4\}$, $E_{2}=\{1,5\}$, $E_{3}=\{1,6\}$, $E_{4}=\{2,5\}$, $E_{5}=\{3,4\}$, $E_{6}=\{4,5\}$, $E_{7}=\{5,7\}$, $E_{8}=\{9\}$, $E_{9}=\{10\}$, $E_{10}=\{11\}$, $E_{11}=\{12\}$.\\ From the circuits it follows also that $\Gamma_{L}$ has 4 connected components which are vertices, namely $\{E_{1}\}$, $\{E_{2}\}$, $\{E_{6}\}$ and $\{E_{10}\}$, 2 connected components which are edges, namely $\{E_{7},E_{9}\}$ and $\{E_{8},E_{11}\}$, and 1 connected component which is a 2-simplex, namely $\{E_{3},E_{4},E_{5}\}$.}
\end{ex1}

The {\em induced subcomplex} $\mathcal{D}'$ of a simplicial complex $\mathcal{D}$ by certain vertices $\mathcal{V}' \subset \mathcal{V}$ is the subcomplex of $\mathcal{D}$ with vertices $\mathcal{V}' $ and $T \subset \mathcal{V}'$ is a simplex of the subcomplex $\mathcal{D}'$ if $T$ is a simplex of $\mathcal{D}$. A subcomplex $H$ of $\mathcal{D}$ is called a {\em spanning subcomplex} if both have exactly the same set of vertices.

Let $F$ be a polynomial in $K[x_{1},\ldots,x_{m}]$. We associate to $F$ the induced subcomplex $\Gamma_{L}(F)$ of $\Gamma_{L}$ consisting of those vertices $E_{i} \in \mathcal{T}_{\rm min}$ with the property: there exists a monomial $M$ in $F$ such that $E_{i}={\rm supp}(M)$.\\ The next theorem provides a necessary condition under which a set of polynomials in the lattice ideal $I_{L}$ generates the radical of $I_{L}$ up to radical.

\begin{thm1} \label{Spanning} Let $K$ be any field. If $rad(I_{L})=rad(F_{1},\ldots,F_{s})$ for some polynomials $F_{1},\ldots,F_{s}$ in $I_{L}$, then $\cup_{i=1}^{s} \Gamma_{L}(F_{i})$ is a spanning subcomplex of $\Gamma_{L}$.
\end{thm1}
\noindent \noindent \textbf{Proof.} Let $E={\rm supp}({\bf x}^{{\bf u}_{+}}) \in \mathcal{T}_{\rm min}$, where $B={\bf x}^{{\bf u}_{+}}-{\bf x}^{{\bf u}_{-}} \in I_{L}$ and the monomial ${\bf x}^{{\bf u}_{+}}$ is indispensable of $I_{L}$. We will prove that there is a monomial $M$ in some $F_{j}$, $1 \leq j \leq s$, such that $E={\rm supp}(M)$. Since $rad(I_{L})=rad(F_{1},\ldots,F_{s})$, there is a power $B^{r}$, $r \geq 1$, that belongs to the ideal $J=(F_{1},\ldots,F_{s}) \subset K[x_{1},\ldots,x_{m}]$. As a consequence there exists a monomial $M$ in some $F_{j}$, $1 \leq j \leq s$, dividing the monomial $({\bf x}^{{\bf u}_{+}})^{r}$, so ${\rm supp}(M) \subseteq {\rm supp}({\bf x}^{{\bf u}_{+}})=E$. Since $F_{j} \in I_{L}$ and $I_{L}$ is generated by binomials, there exists a binomial ${\bf x}^{{\bf v}_{+}}-{\bf x}^{{\bf v}_{-}} \in I_{L}$ such that ${\bf x}^{{\bf v}_{+}}$ divides $M$. But ${\bf x}^{{\bf v}_{+}}$ belongs to $\mathcal{M}_{L}$, so there exists an indispensable monomial $N$ of $I_{L}$ such that $N$ divides ${\bf x}^{{\bf v}_{+}}$. Thus $N$ divides $M$ and therefore ${\rm supp}(N) \subseteq {\rm supp}(M) \subseteq E$. Since $E \in \mathcal{T}_{\rm min}$, we have that $E={\rm supp}(N)$ and therefore $E={\rm supp}(M)$. \hfill $\square$\\

\begin{rem1} \label{RemA} {\rm Let $F$ be an $\mathcal{A}$-homogeneous polynomial of $I_L$, then the simplicial complex $\Gamma_{L}(F)$ is a simplex. To see this suppose that $\Gamma_{L}(F) \neq \emptyset$ and let $T=\{E_{1},\ldots,E_{k}\}$ be the set of vertices of $\Gamma_{L}(F)$. For every $1 \leq i \leq k$ we have that $E_{i} \in \mathcal{T}_{\rm min}$, so, from Theorem \ref{Spanning}, there exists a monomial $M_{i}$, $1 \leq i \leq k$, in $F$ such that $E_{i}={\rm supp}(M_{i})$. But $F$ is $\mathcal{A}$-homogeneous, so the monomials $M_{1},\ldots,M_{k}$ have the same $\mathcal{A}$-degree. By the definition of the simplicial complex $\Gamma_{L}$, we have that $\Gamma_{L}(F)$ is a simplex of $\Gamma_{L}$.}
\end{rem1}

Combining Theorem \ref{Spanning} and Remark \ref{RemA} we take the following corollary.
\begin{cor1} \label{BasicCorol} If $rad(I_{L})=rad(F_{1},\ldots,F_{s})$ for some $\mathcal{A}$-homogeneous polynomials $F_{1},\ldots,F_{s}$ in $I_{L}$, then $\cup_{i=1}^{s} \Gamma_{L}(F_{i})$ is a spanning subcomplex of $\Gamma_{L}$. Furthermore, each $\Gamma_{L}(F_{i})$ is a simplex of $\Gamma_{L}$.
\end{cor1}

\begin{rem1} \label{Remar} {\rm Since any binomial $B=M-N \in I_{L}$ is $\mathcal{A}$-homogeneous, Corollary \ref{BasicCorol} is still valid if we replace every polynomial $F_{i}$, $1 \leq i \leq s$, with a binomial $B_{i}$. Notice that each $\Gamma_{L}(B_{i})$ will be either $1$-simplex, $0$-simplex or the empty set.}
\end{rem1}

Let $\mathcal{D}$ be a simplicial complex on the vertex set $\mathcal{V}$ and let $J$ be a subset of $\Omega=\{0,1,\ldots,{\rm dim}(\mathcal{D})\}$. A set $\mathcal{N}=\{T_{1},\ldots,T_{s}\}$ of simplices of $\mathcal{D}$ is called a $J$-{\em matching} in $\mathcal{D}$ if $T_{k} \cap T_{l}=\emptyset$ for every $1 \leq k,l \leq s$ and ${\rm dim}(T_{k}) \in J$ for every $1 \leq k \leq s$. Let ${\rm supp}(\mathcal{N})=\cup_{i=1}^{s}T_{i}$, which is a subset of the vertices $\mathcal{V}$. A $J$-matching in $\mathcal{D}$ is called a {\em perfect matching} if ${\rm supp}(\mathcal{N})=\mathcal{V}$.\\ A $J$-matching $\mathcal{N}$ in $\mathcal{D}$ is called a {\em maximal} $J$-{\em matching} if ${\rm supp}(\mathcal{N})$ has the maximum possible cardinality among all $J$-matchings.\\ Given a maximal $J$-matching $\mathcal{N}=\{T_{1},\ldots,T_{s}\}$ in $\mathcal{D}$, we shall denote by ${\rm card}(\mathcal{N})$ the cardinality $s$ of the set $\mathcal{N}$. In addition, by $\delta(\mathcal{D})_{J}$ we denote the minimum of the set $$\{{\rm card}(\mathcal{N})| \mathcal{N} \ \textrm{is a maximal} \ J-\textrm{matching in} \ \mathcal{D}.\}$$
It follows from the definitions that if $\mathcal{D}=\cup_{i=1}^{t}\mathcal{D}^{i}$, then $$\delta(\mathcal{D})_{\{0,1\}}=\sum_{i=1}^{t}\delta(\mathcal{D}^{i})_{\{0,1\}}$$ where $\mathcal{D}^{i}$ are the connected components of $\mathcal{D}$.\\ We denote by $c_{\mathcal{D}}$ the smallest number $s$ of simplices $T_{i}$ of $\mathcal{D}$, such that the subcomplex $\bigcup_{i=1}^{s}T_{i}$ is spanning. While by $b_{\mathcal{D}}$ we denote the smallest number $s$ of $1$-simplices or $0$-simplices $T_{i}$ of $\mathcal{D}$, such that the subcomplex $\bigcup_{i=1}^{s}T_{i}$ is spanning.

\begin{thm1} \label{Basic} Let $K$ be any field, then ${\rm bar}(I_L) \geq \delta(\Gamma_{L})_{\{0,1\}}$ and ${\rm ara}_{\mathcal{A}}(I_L) \geq \delta(\Gamma_{L})_{\Omega}$.
\end{thm1}
\noindent \noindent \textbf{Proof.} By Corollary \ref{BasicCorol} and Remark \ref{Remar} we have that ${\rm bar}(I_L) \geq b_{\Gamma_{L}}$ and ${\rm ara}_{\mathcal{A}}(I_L) \geq c_{\Gamma_{L}}$. Now Proposition 3.3 of \cite{KT} asserts that $b_{\Gamma_{{L}}}=\delta(\Gamma_{{L}})_{\{0,1\}}$ and $c_{\Gamma_{{L}}}=\delta(\Gamma_{{L}})_{\Omega}$. Thus ${\rm bar}(I_L) \geq \delta(\Gamma_{{L}})_{\{0,1\}}$ and ${\rm ara}_{\mathcal{A}}(I_L) \geq \delta(\Gamma_{{L}})_{\Omega}$. \hfill $\square$\\

For a vector configuration $\mathcal{B}=\{{\bf b}_{1},\ldots,{\bf b}_{s}\} \subset \mathbb{Z}^{m}$, we denote by $\sigma=pos_{\mathbb{Q}}(\mathcal{B})$ the rational polyhedral cone consisting of all non-negative linear rational combinations of the vectors ${\bf b}_{1},\ldots,{\bf b}_{s}$. Furthermore, $\mathcal{B}$ is called {\em extremal} if for any $\mathcal{S} \subsetneqq \mathcal{B}$ we have $pos_{\mathbb{Q}}(\mathcal{S}) \subsetneqq pos_{\mathbb{Q}}(\mathcal{B})$.\\
 In \cite{KMT} they associated to every lattice ideal $I_{L}$ the rational polyhedral cone $\sigma=pos_{\mathbb{Q}}(\mathcal{A})$ and the simplicial complex $\mathcal{D}_{\sigma}$. As they have shown, combinatorial invariants of $\mathcal{D}_{\sigma}$ provide lower bounds for ${\rm bar}(I_{L})$ and ${\rm ara}_{\mathcal{A}}(I_{L})$. More precisely it holds that ${\rm bar}(I_{L}) \geq \delta(\mathcal{D}_{\sigma})_{\{0,1\}}$ and ${\rm ara}_{\mathcal{A}}(I_{L}) \geq \delta(\mathcal{D}_{\sigma})_{\{\Omega\}}$ (see also Theorem 3.5 of \cite{KT}). Moreover, it was proved in Theorem 4.6 of \cite{KT} that, for an extremal vector configuration $\mathcal{A}$, it holds that $\Delta_{\mathcal{A}}=\mathcal{D}_{\sigma}$, so in this case $\delta(\Gamma_{L})_{\{0,1\}}=\delta(\mathcal{D}_{\sigma})_{\{0,1\}}$ and $\delta(\Gamma_{L})_{\Omega}=\delta(\mathcal{D}_{\sigma})_{\Omega}$. Generally speaking, our lower bounds are essentially different from those derived in \cite{KMT}. The following example shows this fact.

\begin{ex1} {\rm Let $L=ker_{\mathbb{Z}}(\mathcal{A})$ be the lattice of the Example \ref{BasicExample}. We have that $\delta(\Gamma_{L})_{\{0,1\}}=8$, attained by the maximal $\{0,1\}$-matching $$\{\{E_{1}\},\{E_{2}\},\{E_{6}\},\{E_{10}\}, \{E_{7},E_{9}\},\{E_{8},E_{11}\},\{E_{3},E_{4}\},\{E_{5}\}\},$$ so Theorem \ref{Basic} implies that ${\rm bar}(I_{\mathcal{A}}) \geq 8$. Actually ${\rm bar}(I_{\mathcal{A}})=8$, since $\mu(I_{
\mathcal{A}})=8$. Furthermore $\delta(\Gamma_{L})_{\{0,1,2\}}=7$, attained by the maximal $\{0,1,2\}$-matching $$\{\{E_{1}\},\{E_{2}\},\{E_{6}\},\{E_{10}\}, \{E_{7},E_{9}\},\{E_{8},E_{11}\},\{E_{3},E_{4},E_{5}\}\},$$ and therefore, from Theorem \ref{Basic}, the inequality  ${\rm ara}_{\mathcal{A}}(I_{\mathcal{A}}) \geq 7$ holds. Let $\mathcal{Q}$ be the ideal in $K[x_{1},\ldots,x_{9}]$ generated by the binomials $B_{i}$, $1 \leq i \leq 5$, and let $F=x_{1}^{2}x_{6}^2-x_{2}^{2}x_{5}^2+x_{3}^{2}x_{4}^{2}-x_{2}x_{3}x_{6}x_{9} \in \mathcal{Q}$. Then the set of $\mathcal{A}$-homogeneous polynomials $\mathcal{S}=\{F,B_{3},B_{4},B_{5}\} \subset \mathcal{Q}$ generates $rad(\mathcal{Q})$ up to radical, because the polynomials $B_{1}^{3}$ and $B_{2}^{3}$ belong to the ideal generated by the polynomials in $\mathcal{S}$. Consequently $I_{\mathcal{A}}$ is generated up to radical by seven $\mathcal{A}$-homogenous polynomials, namely $B_{i}$, $3 \leq i \leq 8$, and $F$. Thus ${\rm ara}_{\mathcal{A}}(I_{\mathcal{A}})=7$.\\ Notice that $\mathcal{A}$ is not an extremal configuration. Actually $\mathcal{B}=\{{\bf a}_{1},\ldots,{\bf a}_{8},{\bf a}_{9}\} \subsetneqq \mathcal{A}$ is an extremal vector configuration. To compute the simplicial complex $\Delta_{\mathcal{B}}=\mathcal{D}_{\sigma}$ one should find the circuits of the toric ideal $I_{\mathcal{B}}$. Proposition 4.13 in \cite{St} asserts that the circuits of $I_{\mathcal{B}}$ are $\mathcal{C}_{\mathcal{B}}=\mathcal{C}_{\mathcal{A}} \cap K[x_{1},\ldots,x_{9}]$. The simplicial complex $\Gamma_{ker_{\mathbb{Z}}(\mathcal{B})}$ has 9 vertices, namely $E_{i}$, $1 \leq i \leq 6$, $E_{7}'=\{2,9\}$, $E_{8}'=\{3,9\}$ and $E_{9}'=\{6,9\}$. Furthermore it has 3 connected components which are edges, namely $\{E_{1},E_{7}'\}$, $\{E_{2},E_{8}'\}$ and $\{E_{6},E_{9}'\}$, and also 1 connected component which is a 2-simplex, namely $\{E_{3},E_{4},E_{5}\}$. It follows easily that $\delta(\mathcal{D}_{\sigma})_{\{0,1\}}=\delta(\Delta_{\mathcal{B}})_{\{0,1\}}=5$ and also $\delta(\mathcal{D}_{\sigma})_{\{0,1,2\}}=\delta(\Delta_{\mathcal{B}})_{\{0,1\}}=4$. Remark that ${\rm ht}(I_{\mathcal{A}})=6$.}
\end{ex1}

\begin{prop1} \label{Integer} Let $q$ be the number of vertices of $\mathcal{T}_{\rm min}$, then ${\rm bar}(I_{L}) \geq \lceil \frac{q}{2} \rceil$.
\end{prop1}
\noindent \textbf{Proof.} By Remark 2.5 in \cite{KT} every maximal $\{0,1\}$-matching in $\Gamma_{L}$ is perfect. Clearly $\delta(\Gamma_{L})_{\{0,1\}} \geq \lceil \frac{q}{2} \rceil$ and therefore we have, from Theorem \ref{Basic}, that ${\rm bar}(I_{L}) \geq \lceil \frac{q}{2} \rceil$. \hfill $\square$\\

In this work our basic aim is to study when the equality ${\rm bar}(I_{L})=\mu(I_{L})$ holds. Of particular interest is the case that $I_{L}$ has a generating set $\{B_{1},\ldots,B_{t}\}$ such that every binomial $B_{i}$ is a difference of squarefree monomials. The next theorem asserts that the above equality holds for such ideals, under the assumption that the lattice ideal $I_{L}$ is generated by its indispensable.

\begin{thm1} Suppose that the lattice ideal $I_{L}$ has a binomial generating set $\{B_{1},\ldots,B_{t}\}$ such that every $B_{i}$, $1 \leq i \leq t$, is a difference of squarefree monomials. If $I_{L}$ has a unique minimal system of binomial generators, then ${\rm bar}(I_{L})=\mu(I_{L})$.
\end{thm1}
\noindent \textbf{Proof.}  Since $I_{L}$ is generated by binomials which are differences of squarefree monomials, every indispensable monomial of $I_{L}$ is squarefree. First we prove that the support of every indispensable monomial $M$ of $I_{L}$ belongs to $\mathcal{T}_{\rm min}$. If there exists an indispensable monomial $N$ of $I_{L}$ such that ${\rm supp}(N) \subsetneqq {\rm supp}(M)$, then $N$ divides $M$ and $N \neq M$, a contradiction to the fact that $M$ is indispensable. Let $\mathcal{P}$ be the unique minimal binomial generating set of $I_{L}$. We claim that for an indispensable monomial $M$ of $I_{L}$ there exists exactly one binomial $B \in \mathcal{P}$ such that $M$ is a monomial of $B$. Let $B=M-N$ and suppose that there exists another binomial $B' \in \mathcal{P}$ such that $B'=M-N'$. Then we can replace $B$ by $B'$ and $N-N'$ in $\mathcal{P}$, thus obtaining a system of generators of $I_{L}$ not containing $B$ which is not possible by Definition \ref{Indisp}. Let $q$ be the number of vertices of $\mathcal{T}_{\rm min}$ and $s=\mu(I_{L})$, then $q=2s$ and therefore we have, from Proposition \ref{Integer}, that ${\rm bar}(I_{L}) \geq s$. Consequently ${\rm bar}(I_{L})=s$. \hfill $\square$\\



\section{The case of toric ideals associated with graphs}

In this section we consider a special class of lattice ideals, namely toric ideals associated with graphs. In the sequel, all graphs under consideration are finite, simple and connected. Recall that a simple graph is an abstract simplicial complex consisting only of vertices and edges. To every graph $G$ is associated the toric ideal $I_{\mathcal{A}_{G}}$. We study the equality ${\rm bar}(I_{\mathcal{A}_{G}})=\mu(I_{\mathcal{A}_{G}})$, when $I_{\mathcal{A}_{G}}$ has a generating set $\{B_{1},\ldots,B_{t}\}$ such that every binomial $B_{i}$ is a difference of squarefree monomials.\\


\subsection{Basics on toric ideals of graphs}\mbox{}\par

\medskip

Let $G$ be a graph on the vertex set $\mathcal{V}(G)=\{v_{1},\ldots,v_{n}\}$ with edges $\mathcal{E}(G)=\{e_{1},\ldots,e_{m}\}$. Consider one variable $x_{i}$ for each $e_{i}$ and form the polynomial ring $K[x_{1},\ldots,x_{m}]$ over any field $K$. To every edge $e=\{v_{i},v_{j}\} \in \mathcal{E}(G)$ we associate the vector ${\bf a}_{e}=(0,\ldots,0,1,0,\ldots,0,1,0,\ldots,0) \in \mathbb{Z}^{n}$ with exactly two 1's, which are in $i$ and $j$ position, and the rest of its entries equal to zero. Let $\mathcal{A}_{G}=\{{\bf a}_{e}|e \in \mathcal{E}(G)\}$ and consider the toric ideal $I_{\mathcal{A}_{G}} \subset K[x_{1},\ldots,x_{m}]$.

\begin{not1} {\rm For the sake of simplicity we are going to write ${\bf a}_{i}$, $1 \leq i \leq m$, instead of ${\bf a}_{e_{i}}$ and $\Gamma_{G}$ for the simplicial complex $\Gamma_{ker_{\mathbb{Z}}(\mathcal{A}_{G})}$.}
\end{not1}

A {\em walk} of length $s$ of $G$ is a finite sequence of the form $$w=(\{v_{1},v_{2}\},\{v_{2},v_{3}\},\ldots,\{v_{s},v_{s+1}\}).$$ We say that the walk is {\em closed} if $v_{1}=v_{s+1}$. An {\em even} (respectively {\em odd}) closed walk is a closed walk of even (respectively odd) length.\ A {\em cycle} of $G$ is a closed walk $w=(\{v_{1},v_{2}\},\{v_{2},v_{3}\},\ldots,\{v_{s},v_{1}\})$ with $v_{i} \neq v_{j}$, for every $1 \leq i<j \leq s$.\\ For an even closed walk $w=(e_{i_1},e_{i_2},\ldots,e_{i_{2s}})$ of $G$ with each $e_{k} \in \mathcal{E}({G})$, it holds that $$\phi(\prod_{k=1}^{s} x_{i_{2k-1}})=\phi(\prod_{k=1}^{s} x_{i_{2k}})$$ and therefore the binomial $$B_{w}:=\prod_{k=1}^{s} x_{i_{2k-1}}-\prod_{k=1}^{s} x_{i_{2k}}$$ belongs to $I_{\mathcal{A}_{G}}$. We often employ the abbreviated notation $$B_{w}=B_{w}^{(+)}-B_{w}^{(-)},$$ where $$B_{w}^{(+)}=\prod_{k=1}^{s} x_{i_{2k-1}}, \ B_{w}^{(-)}=\prod_{k=1}^{s} x_{i_{2k}}.$$From Proposition 3.1 in \cite{Vil} we have that every toric ideal $I_{\mathcal{A}_{G}}$ is generated by binomials of the above form.
\begin{rem1} \label{Extremal} {\rm The toric ideal $I_{\mathcal{A}_{G}}$ has no binomials of the form $B=x_{i}^{u_{i}}-{\bf x}^{\bf v}$ where $i \notin {\rm supp}({\bf x}^{\bf v})$. If $I_{\mathcal{A}_{G}}$ has such a binomial $B$, then ${\rm deg}_{\mathcal{A}_{G}}(x_{i}^{u_{i}})={\rm deg}_{\mathcal{A}_{G}}({\bf x}^{\bf v})$. Combining this fact together with that the entries of every ${\bf a}_{j}$, $1 \leq j \leq m$, are either $0$ or $1$ and exactly two of them are equal to $1$, we arrive at a contradiction.}
\end{rem1}

A binomial $B \in I_{\mathcal{A}_{G}}$ is called {\em minimal} if it belongs to a minimal system of binomial generators of $I_{\mathcal{A}_{G}}$. Every minimal binomial is primitive, see \cite{MS}, \cite{St}. Recall that an irreducible binomial ${\bf x}^{{\bf u}_+}-{\bf x}^{{\bf u}_-} \in I_{\mathcal{A}_{G}}$ is called {\em primitive} if there exists no other binomial ${\bf x}^{{\bf v}_+}-{\bf x}^{{\bf v}_-} \in I_{\mathcal{A}_{G}}$ such that ${\bf x}^{{\bf v}_+}$ divides ${\bf x}^{{\bf u}_+}$ and ${\bf x}^{{\bf v}_-}$ divides ${\bf x}^{{\bf u}_-}$. For a primitive binomial $B={\bf x}^{{\bf u}_+}-{\bf x}^{{\bf u}_-} \in I_{\mathcal{A}_{G}}$ we have, from Lemma 3.2 in \cite{OH1}, that $B=B_{w}$ for an even closed walk $w$ of certain type. An even closed walk $w=(e_{i_1},\ldots,e_{i_{2s}})$ of $G$ is called {\em primitive} if there exists no even closed walk of $G$ of the form $(e_{j_1},\ldots,e_{j_{2t}})$ with $1 \leq t<s$ such that each $j_{2k-1}$ belongs to $\{i_{1},i_{3},\ldots,i_{2s-1}\}$, each $j_{2k}$ belongs to $\{i_{2},i_{4},\ldots,i_{2s}\}$ and $j_{2k-1} \neq j_{l}$ for all $1 \leq k \leq t$ and for all $1 \leq l \leq t$. The walk $w$ is primitive if and only if the binomial $B_{w}$ is primitive.

Every circuit $B \in I_{\mathcal{A}_{G}}$ is also a primitive binomial, so $B=B_{w}$ for an even closed walk $w$ of $G$. The next theorem provides a characterization of all even closed walks $w$ such that $B_{w}$ is a circuit.

\begin{thm1} \label{Circuits} (\cite{Vil}) Let $G$ be a graph. Then a binomial $B \in I_{\mathcal{A}_{G}}$ is a circuit if and only if $B=B_{w}$ where \begin{enumerate} \item $w$ is an even cycle or \item two odd cycles intersecting in exactly one vertex or \item two vertex disjoint odd cycles joined by a path.
\end{enumerate}
\end{thm1}
\begin{rem1} {\rm Let $B_{w}$ be a circuit. Then the monomials $B_{w}^{(+)}$, $B_{w}^{(-)}$ are squarefree if and only if $w$ is an even cycle or two odd cycles intersecting in exactly one vertex.}
\end{rem1}

For the rest of this section we recall some fundamental material from \cite{RTT}. A {\em cut vertex} $v$ in a graph $G$ is a vertex, such that if $v$ is removed, the number of connected components of $G$ increases. A connected graph is said to be {\em biconnected} if it does not contain a cut vertex. A maximal biconnected subgraph of a graph is called a {\em block}.\\ Using the fact that every primitive binomial $B_{w}$ is irreducible, we deduce that the set of edges of $w=(e_{1},\ldots,e_{2s})$ has a partition into two sets, namely $w^{+}=\{e_{1},e_{3},\ldots,e_{2s-1}\}$ and $w^{-}=\{e_{2},e_{4},\ldots,e_{2s}\}$. The edges of $w^{+}$ are called {\em odd edges} of $w$, while the edges of $w^{-}$ are called {\em even} edges of $w$.

Given a primitive walk $$w=(e_{1}=\{v_{1},v_{2}\},e_{2}=\{v_{2},v_{3}\},\ldots,e_{2s}=\{v_{2s},v_{1}\})$$ of $G$ which has a chord $e=\{v_{k},v_{l}\}$ with $1 \leq k<l \leq 2s$, we have that $e$ breaks $w$ in the walks $\gamma_{1}=(e_{1},\ldots,e_{k-1},e,e_{l},\ldots,e_{2s})$ and $\gamma_{2}=(e_{k},\ldots,e_{l-1},e)$. The chord $e$ is called {\em bridge} of $w$ if there are two different blocks $\mathcal{B}_{1}$, $\mathcal{B}_{2}$ of $w$ such that $v_{k} \in \mathcal{B}_{1}$ and $v_{l} \in \mathcal{B}_{2}$. Furthermore, $e$ is called {\em odd} if it is not a bridge and both $\gamma_{1}$, $\gamma_{2}$ are odd walks. Notice that if $e$ is an odd chord of $w$, then $l-k$ is even.

\begin{def1} \label{Cross} {\rm Let $w=(\{v_{i_1},v_{i_2}\},\{v_{i_2},v_{i_3}\},\ldots,\{v_{i_{2s}},v_{i_1}\})$ be a primitive walk. Given two odd chords $e_{1}=\{v_{i_t},v_{i_j}\}$ and $e_{2}=\{v_{i_{t'}},v_{i_{j'}}\}$ with $1 \leq t<j \leq 2s$ and $1 \leq t'<j' \leq 2s$, we say that \begin{enumerate} \item $e_{1}$ and $e_{2}$ cross effectively in $w$ if $t'-t$ is odd and either $t<t'<j<j'$ or $t'<t<j'<j$. \item $e_{1}$ and $e_{2}$ cross strongly effectively in $w$ if they cross effectively and they don't form an $\mathcal{F}_{4}$ in $w$.
\end{enumerate}}
\end{def1}

\begin{def1} {\rm Let $w$ be a primitive walk of $G$. We call an $\mathcal{F}_{4}$ of $w$ an even cycle $\xi=(e_{i},e_{j},e_{k},e_{l})$ of length 4 consisting of two edges $e_{i}$, $e_{k}$ of $w$, which are both even or both odd, and the odd chords $e_{j}$ and $e_{l}$ which cross effectively in $w$.}
\end{def1}
 \begin{rem1} \label{F4} {\rm (1) If $\xi=(e_{i},e_{j},e_{k},e_{l})$ is an $\mathcal{F}_{4}$ of a primitive walk $w$, where $e_{j}$ and $e_{l}$ are two odd chords which cross effectively in $w$, then $x_{i}x_{k}$ divides exactly one of the monomials $B_{w}^{(+)}$ and $B_{w}^{(-)}$.\\
 (2) If $B_{w}$ is a minimal binomial which is not indispensable, then combining Theorem 4.13, Proposition 4.10 and Theorem 4.14 in \cite{RTT} we deduce that the walk $w$ has at least one $\mathcal{F}_{4}$.}
\end{rem1}

Let $\xi=(e_{i},e_{j},e_{k},e_{l})$ be an $\mathcal{F}_{4}$ of a primitive walk $w$, where $e_{j}$ and $e_{l}$ are odd chords of $w$. The walk $w$ can be written as $w=(w_{1},e_{i},w_{2},e_{k})$, where $w_{1}$, $w_{2}$ are walks in $G$. Notice that every walk $\gamma$ can be can be regarded as a subgraph of $G$ with vertices the vertices of the walk and edges the edges of the walk $\gamma$. The $\mathcal{F}_{4}$ induces a partition of the vertices of $w$ into the sets $\mathcal{V}(w_1)$, $\mathcal{V}(w_{2})$. Where $\mathcal{V}(w_1)$, $\mathcal{V}(w_{2})$ denote the set of vertices of $w_{1}$ and $w_{2}$, respectively. We say that an odd chord $e$ of the primitive walk $w$ {\em crosses} the $\mathcal{F}_{4}$ if one of the vertices of $e$ belongs to $\mathcal{V}(w_{1})$, the other belongs to $\mathcal{V}(w_{2})$ and $e$ is different from $e_{j}$, $e_{l}$.

\begin{ex1} {\rm Let $G$ be the graph on the vertex set $\{v_{1},\ldots,v_{8}\}$ with edges $e_{i}=\{v_{i},v_{i+1}\}$, $1 \leq i \leq 7$, $e_{8}=\{v_{1},v_{8}\}$, $e_{9}=\{v_{1},v_{5}\}$, $e_{10}=\{v_{2},v_{4}\}$, $e_{11}=\{v_{4},v_{6}\}$ and $e_{12}=\{v_{5},v_{7}\}$. Consider the even cycle $w=(e_{1},\ldots,e_{8})$ which has four odd chords, namely $e_{9}$, $e_{10}$, $e_{11}$ and $e_{12}$. For instance the odd chords $e_{9}$ and $e_{10}$ don't cross effectively. On the contrary, the odd chords $e_{11}$, $e_{12}$ cross effectively and they form an $\mathcal{F}_{4}$ of $w$, namely the even cycle $\xi=(e_{4},e_{12},e_{6},e_{11})$. The even cycle $w$ can be written as $w=(w_{1},e_{4},w_{2},e_{6})$, where $w_{1}=(e_{7},e_{8},e_{1},e_{2},e_{3})$ and $w_{2}=e_{5}$. Now the odd chord $e_{9}$ crosses the $\mathcal{F}_{4}$, since $v_{1} \in \mathcal{V}(w_{1})$ and $v_{5} \in \mathcal{V}(w_{2})$.}
\end{ex1}


\subsection{Binomial arithmetical rank of the toric ideal associated with a graph} \mbox{}\par

\medskip

Recently H. Ohsugi and T. Hibi (\cite{OH}) provided a characterization of all graphs $G$ such that the toric ideal $I_{\mathcal{A}_{G}}$ is generated by circuits of the form ${\bf x}^{{\bf u}_{+}}-{\bf x}^{{\bf u}_{-}}$, where both monomials ${\bf x}^{{\bf u}_{+}}$ and ${\bf x}^{{\bf u}_{-}}$ are squarefree. More precisely they proved that the following are equivalent:
\begin{enumerate} \item $I_{\mathcal{A}_{G}}$ is generated by circuits of the form ${\bf x}^{{\bf u}_{+}}-{\bf x}^{{\bf u}_{-}}$, where both monomials ${\bf x}^{{\bf u}_{+}}$ and ${\bf x}^{{\bf u}_{-}}$ are squarefree. \item There is no induced subgraph of $G$ consisting of two odd cycles vertex disjoint joined by a path of length $\geq 1$.
\end{enumerate}
From now on every graph $G$, unless otherwise stated, will satisfy the condition $(\sharp)$: {\em There is no induced subgraph of $G$ consisting of two odd cycles vertex disjoint joined by a path of length $\geq 1$.}

\begin{ex1} \label{Wheel} {\rm Let $\xi_{n}=(\{v_{1},v_{2}\},\{v_{2},v_{3}\},\ldots,\{v_{n-1},v_{n}\},\{v_{n},v_{1}\})$ be a cycle of length $n \geq 3$. The wheel graph $W_{n+1}$ on the vertex set $\{v_{1},\ldots,v_{n},v_{n+1}\}$ is the graph with edges all the edges of $\xi_{n}$ and also $\{v_{i},v_{n+1}\}$ is an edge of $W_{n+1}$, for every $1 \leq i \leq n$. If $n$ is even, then $v_{n+1}$ is a vertex of every odd cycle of $W_{n+1}$. If $n$ is odd, then any odd cycle of $W_{n+1}$ either coincides with $\xi_{n}$ or has at least 3 vertices, namely $v_{n+1}$ and 2 vertices of $\xi_{n}$. In both cases $W_{n+1}$ has no two odd cycles vertex disjoint, so $W_{n+1}$ satisfies $(\sharp)$.}
\end{ex1}

Recall that the vertices of the simplicial complex $\Gamma_{G}$ are exactly the elements of $\mathcal{T}_{\rm min}$. By Remark \ref{RemCircuE} (2) there is an edge $\{E_{i},E_{j}\}$ of $\Gamma_{G}$ if and only if there exists a circuit ${\bf x}^{{\bf u}_{+}}-{\bf x}^{{\bf u}_{-}} \in I_{\mathcal{A}_{G}}$ such that $E_{i}={\rm supp}({\bf x}^{{\bf u}_{+}})$ and $E_{j}={\rm supp}({\bf x}^{{\bf u}_{-}})$. We will detect the structure of every connected component of the simplicial complex $\Gamma_{G}$.

\begin{prop1} \label{Edge} Let $E={\rm supp}(M_{i})$ and $E'={\rm supp}(M_{j})$ be two vertices of $\Gamma_{G}$, where $M_{i}$, $M_{j}$ are indispensable monomials of $I_{\mathcal{A}_{G}}$. Then \begin{enumerate} \item[(1)] $\{E,E'\}$ is an edge of $\Gamma_{G}$ if and only if there exists a circuit $B_{w}=M_{i}-M_{j} \in I_{\mathcal{A}_{G}}$. \item[(2)] An edge $\{E,E'\}$ is a connected component of $\Gamma_{G}$ if and only if there is an indispensable binomial $B_{w}=M_{i}-M_{j} \in I_{\mathcal{A}_{G}}$ with $E={\rm supp}(M_{i})$, $E'={\rm supp}(M_{j})$.
\end{enumerate}
\end{prop1}
\noindent \noindent \textbf{Proof.} (1) ($\Leftarrow$) If there exists a circuit $B_{w}=M_{i}-M_{j} \in I_{\mathcal{A}_{G}}$, then we have, from the definition of the complex $\Gamma_{G}$, that $\{E,E'\}$ is an edge of $\Gamma_{G}$.\\
($\Rightarrow$) Assume that $\{E,E'\}$ is an edge of $\Gamma_{G}$, then there exists a circuit $B_{w}=B_{w}^{(+)}-B_{w}^{(-)} \in I_{\mathcal{A}_{G}}$ with ${\rm supp}(B_{w}^{(+)})=E \in \mathcal{T}_{\rm min}$ and ${\rm supp}(B_{w}^{(-)})=E' \in \mathcal{T}_{\rm min}$. It is enough to prove that the monomials $B_{w}^{(+)}$, $B_{w}^{(-)}$ are indispensable of $I_{\mathcal{A}_{G}}$. If $w$ is an even cycle or two odd cycles intersecting in exactly one vertex, then both monomials $B_{w}^{(+)}$ and $B_{w}^{(-)}$ are squarefree and therefore they are indispensable, since ${\rm supp}(B_{w}^{(+)})$, ${\rm supp}(B_{w}^{(-)}) \in \mathcal{T}_{\rm min}$. Thus necessarily in this case $B_{w}^{(+)}=M_{i}$ and $B_{w}^{(-)}=M_{j}$. Let us now assume that $w$ consists of two vertex disjoint odd cycles $\xi_{1}$, $\xi_{2}$ joined by a path $\gamma_{1}=(e_{1},\ldots,e_{r})$ of length $r \geq 1$ connecting one vertex $i$ of $\xi_{1}$ with one vertex $j$ of $\xi_{2}$. We distinguish the following cases: \begin{enumerate} \item[(i)] $r=1$, so $e_{r}=\{i,j\}$. Since $G$ satisfies condition $(\sharp)$, there is an edge $e=\{p,q\} \neq e_{r}$ (i.e. $p \neq i$ or/and $q \neq j$) between one vertex $p$ of $\xi_{1}$ and one vertex $q$ of $\xi_{2}$. Let, say, that $p \neq i$, then there are two paths in $\xi_{1}$ joining $p$ with $i$. Denote by $V_{1}$, $V_{2}$ the paths of even and odd length, respectively, joining $p$ with $i$. In case that $q=j$ we consider the even cycle $\gamma=(p,V_{1},i,e_{r},j,e,p)$ of $G$. Notice that $B_{\gamma}$ is a circuit. Without loss of generality we can assume that $w=(p,V_{1},i,e_{r},j,\xi_{2},j,e_{r},i,V_{2},p)$, then ${\rm supp}(B_{\gamma}^{(+)}) \subsetneqq {\rm supp}(B_{w}^{(+)})$ and therefore, from Lemma \ref{Nosubset}, it holds that $E \notin \mathcal{T}_{\rm min}$, a contradiction. Assume, now, that $q \neq j$. Let $W_{1}$, $W_{2}$ be paths in $\xi_2$ of even and odd length, respectively, joining $q$ with $j$. Consider the even cycle $\gamma=(p,V_{1},i,e_{r},j,W_{1},q,e,p)$. Without loss of generality we can assume that $w=(p,V_{1},i,e_{r},j,W_{1},q,W_{2},j,e_{r},i,V_{2},p)$. Then ${\rm supp}(B_{\gamma}^{(+)}) \subsetneqq {\rm supp}(B_{w}^{(+)})$, a contradiction. \item[(ii)] $r>1$. Suppose first that there exists an edge of $G$ joining a vertex of $\xi_{1}$ with a vertex of $\xi_{2}$. Since $G$ satisfies condition $(\sharp)$, there is no induced subgraph of $G$ consisting of two odd cycles vertex disjoint joined by an edge. Thus there exists at least one edge $e=\{p,q\}$ joining $\xi_{1}$ and $\xi_{2}$, where $p \neq i$ or/and $q \neq j$. Let, say, that $p \neq i$ and assume that $q \neq j$. Let $V_{1}$, $V_{2}$ be paths in $\xi_{1}$ of even and odd length, respectively, joining $p$ with $i$. Let $W_{1}$, $W_{2}$ be paths in $\xi_{2}$ of even and odd length, respectively, joining $q$ with $j$. If the length of $\gamma_{1}$ is odd, we consider the even cycle $\gamma=(p,V_{2},i,\gamma_{1},j,W_{2},q,e,p)$. Assuming that $w=(p,V_{2},i,\gamma_{1},j,W_{2},q,W_{1},j,\gamma_{2},i,V_{1},p)$, where $\gamma_{2}=(e_{r},\ldots,e_{1})$, we have that ${\rm supp}(B_{\gamma}^{(+)}) \subsetneqq {\rm supp}(B_{w}^{(+)})$, a contradiction. If the length of $\gamma_{1}$ is even, we consider the even cycle $\gamma=(p,V_{1},i,\gamma_{1},j,W_{2},q,e,p)$. Assuming that $$w=(p,V_{1},i,\gamma_{1},j,W_{2},q,W_{1},j,\gamma_{2},i,V_{2},p),$$ we have ${\rm supp}(B_{\gamma}^{(+)}) \subsetneqq {\rm supp}(B_{w}^{(+)})$, a contradiction. Using similar arguments we can arrive at a contradiction when $q=j$.\\ Suppose, now, that there exists no such edge. Then there exists an edge of $G$ joining a vertex $p$ of $\xi_{1}$ with a vertex $q$ ($\neq i$) of $\gamma_{1}=(e_{1},\ldots,e_{r})$ and $e=\{p,q\}$ does not belong to $w$. Let $V_{1}$ be a path in $\xi_{1}$ joining $p$ with $i$ and $W_{1}$ be a path in $\xi_{2}$ joining $i$ with $q$. Without loss of generality we can assume that the path $(V_{1},W_{1})$ is odd. Consider the even cycle $\gamma=(p,V_{1},i,W_{1},q,e,p)$. Notice that $B_{\gamma}$ is a circuit. Then ${\rm supp}(B_{\gamma}^{(+)}) \subsetneqq {\rm supp}(B_{w}^{(+)})$ or ${\rm supp}(B_{\gamma}^{(+)}) \subsetneqq {\rm supp}(B_{w}^{(-)})$ and therefore, from Lemma \ref{Nosubset}, it holds that $E \notin \mathcal{T}_{\rm min}$ or $E' \notin \mathcal{T}_{\rm min}$, a contradiction.
\end{enumerate}
(2) ($\Leftarrow$)  Suppose that the edge $\{E,E'\}$ is not a connected component of $\Gamma_{G}$ and let $E''={\rm supp}(M_{k})$ such that $\{E',E''\}$ is an edge of $\Gamma_{G}$. Then the binomials $B_{w}, M_{j}-M_{k}$ and $M_{i}-M_{k}$ belong to $I_{\mathcal{A}_{G}}$ and therefore all the monomials $M_{i}$, $M_{j}$ and $M_{k}$ have the same $\mathcal{A}_{G}$-degree. Thus $\{M_{i},M_{j},M_{k}\}$ is a face of the indispensable complex $\Delta_{{\rm ind}(\mathcal{A}_{G})}$. But the binomial $B_{w}$ is indispensable of $I_{\mathcal{A}_{G}}$, so we have, from Theorem 3.4 in \cite{CKT}, that $\{M_{i},M_{j}\}$ is a facet of $\Delta_{{\rm ind}(\mathcal{A}_{G})}$ and therefore $\{M_{i},M_{j},M_{k}\}$ can't be a face of $\Delta_{{\rm ind}(\mathcal{A}_{G})}$. Consequently, the edge $\{E,E'\}$ is a connected component of $\Gamma_{G}$.\\
($\Rightarrow$) Suppose that the edge $\{E,E'\}$ is a connected component of $\Gamma_{G}$. Then there is a circuit $B_{w}=M_{i}-M_{j} \in I_{\mathcal{A}_{G}}$. Since $M_{i}$, $M_{j}$ are indispensable monomials, we have that $B_{w}$ is a minimal binomial of $I_{\mathcal{A}_{G}}$, see Theorem 1.8 of \cite{KO}. If $B_{w}$ is not indispensable, then we have, from Remark \ref{F4} (2), that the walk $w$ has at least one $\mathcal{F}_{4}$, namely an even cycle $\xi=(e_{1},e_{2},e_{3},e_{4})$ where $e_{2}$ and $e_{4}$ are odd chords of $w$. So the circuit $B_{\xi}=x_{1}x_{3}-x_{2}x_{4}$ belongs to $I_{\mathcal{A}_{G}}$ and also the monomial $x_{1}x_{3}$ divides one of the monomials $M_{i}$ and $M_{j}$, say $M_{i}$. Since the monomial $M_{i}$ is indispensable of $I_{\mathcal{A}_G}$, we have that $M_{i}=x_{1}x_{3}$. Thus $M_{j}$ is quadratic and also, from Remark \ref{Extremal}, the support of the monomial $N=x_{2}x_{4}$ belongs to $\mathcal{T}_{\rm min}$. Now $\{E,E',E''={\rm supp}(N)\}$ is a 2-simplex of $\Gamma_{G}$, a contradiction to the fact that $\{E,E'\}$ is a connected component of $\Gamma_{G}$. \hfill $\square$\\

\begin{thm1} (1) Let $M$ be an indispensable monomial of $I_{\mathcal{A}_{G}}$ that is not quadratic. Then $\{{\rm supp}(M)\}$ is a connected component of $\Gamma_{G}$ if and only if every walk $w$, such that $M$ is a monomial of $B_{w}$, has an $\mathcal{F}_{4}$.\\(2) Every connected component of $\Gamma_{G}$ is either a vertex, an edge or a 2-simplex.
\end{thm1}
\noindent \noindent \textbf{Proof.} (1) Let $E={\rm supp}(M)$. Suppose first that every walk $w$, such that $B_{w}^{(+)}=M$, has an $\mathcal{F}_{4}$. Let us assume that there exists an edge $\{E,E'\}$ of $\Gamma_{G}$, where $E'={\rm supp}(N) \in \mathcal{T}_{\rm min}$ and $N$ is an indispensable monomial of $I_{\mathcal{A}_{G}}$. Then we have, from Proposition \ref{Edge} (1), that the binomial $M-N \in I_{\mathcal{A}_{G}}$ is a circuit, so it is of the form $B_{\gamma}$ for an even closed walk $\gamma$. Notice that the monomials $M$, $N$ are squarefree and both of them they are not quadratic. From the assumption $\gamma$ has an $\mathcal{F}_{4}$, namely $\xi=(e_{1},e_{2},e_{3},e_{4})$ where $e_{2}$ and $e_{4}$ are odd chords of $\gamma$. So the binomial $B_{\xi}=x_{1}x_{3}-x_{2}x_{4} \in I_{\mathcal{A}_{G}}$ is a circuit and also $x_{1}x_{3}$ divides one of the monomials $M$ and $N$, a contradiction to the fact that $M$, $N$ are non-quadratic indispensable monomials. Conversely assume that $\{E\}$ is a connected component of $\Gamma_{G}$. Let $w$ be an even closed walk such that $M$ is a monomial of $B_{w}$, i.e. $B_{w}=M-{\bf x}^{\bf v}$. Then $M$ is indispensable of $I_{\mathcal{A}_G}$ and therefore, from Theorem 1.8 of \cite{KO}, $B_{w}$ is a minimal binomial of $I_{\mathcal{A}_{G}}$. If ${\bf x}^{\bf v}$ is indispensable, then it is squarefree and therefore ${\rm supp}({\bf x}^{\bf v}) \in \mathcal{T}_{\rm min}$. Thus $\{E,{\rm supp}({\bf x}^{\bf v})\}$ is an edge of $\Gamma_{G}$, a contradiction to the fact that $\{E\}$ is a connected component of $\Gamma_{G}$. Consequently ${\bf x}^{\bf v}$ is not indispensable, so $B_{w}$ is not an indispensable binomial of $I_{\mathcal{A}_{G}}$. By Remark \ref{F4} (2) the walk $w$ has at least one $\mathcal{F}_{4}$.\\
(2) First we will show that $\{E,E',E''\}$ is a 2-simplex of $\Gamma_{G}$ if and only if there are quadratic binomials $M_{i}-M_{j}$, $M_{j}-M_{k}$, $M_{i}-M_{k}$ in $I_{\mathcal{A}_{G}}$ with ${\rm supp}(M_{i})=E$, ${\rm supp}(M_{j})=E'$ and ${\rm supp}(M_{k})=E''$. The if implication is easily derived form the fact that ${\rm deg}_{\mathcal{A}_{G}}(M_{i})={\rm deg}_{\mathcal{A}_{G}}(M_{j})={\rm deg}_{\mathcal{A}_{G}}(M_{k})$. Conversely assume that $\{E,E',E''\}$ is a 2-simplex of $\Gamma_{G}$. So there exist indispensable monomials $M_{i}$, $M_{j}$, $M_{k}$ of $I_{\mathcal{A}_{G}}$ with $E={\rm supp}(M_{i})$, $E'={\rm supp}(M_{j})$, $E''={\rm supp}(M_{k})$ such that all binomials $M_{i}-M_{j} \in I_{\mathcal{A}_{G}}$, $M_{j}-M_{k} \in I_{\mathcal{A}_{G}}$ and $M_{i}-M_{k} \in I_{\mathcal{A}_{G}}$ are circuits. Since $M_{i}-M_{j}=B_{w}$ is a minimal binomial of $I_{\mathcal{A}_{G}}$ which is not indispensable, the walk $w$ has an $\mathcal{F}_{4}$ and therefore the indispensable monomials $M_{i}$, $M_{j}$ are quadratic, as well as the monomial $M_{k}$.\\ If for instance there exists an $E'''={\rm supp}(M_{l}) \in \mathcal{T}_{\rm min}$ such that $\{E,E'''\}$ is an edge, then the binomial $M_{i}-M_{l} \in I_{\mathcal{A}_{G}}$ is a circuit and also the monomial $M_{l}$ is quadratic, since $M_{i}$ is quadratic. Thus $M_{l}$ equals either $M_{j}$ or $M_{k}$, see the proof of Proposition 3.4 (2) in \cite{K}. Consequently, every connected component of $\Gamma_{G}$ is either a vertex, an edge or a 2-simplex.  \hfill $\square$\\

The following example demonstrates that there are graphs $G$, such that $\Gamma_{G}$ has a connected component which is a vertex.

\begin{ex1} {\rm Let $G$ be the graph on the vertex set $\{v_{1},\ldots,v_{6}\}$ with edges $e_{i}=\{v_{i},v_{i+1}\}$, $1 \leq i \leq 5$, $e_{6}=\{v_{1},v_{6}\}$, $e_{7}=\{v_{1},v_{3}\}$, $e_{8}=\{v_{2},v_{4}\}$. The circuits are $B_{w_1}=x_{1}x_{3}-x_{7}x_{8}$, $B_{w_2}=x_{2}x_{4}x_{6}-x_{5}x_{7}x_{8}$ and $B_{w_3}=x_{1}x_{3}x_{5}-x_{2}x_{4}x_{6}$. Actually the toric ideal $I_{\mathcal{A}_{G}}$ is minimally generated by the binomials $B_{w_1}$ and $B_{w_2}$. Thus $$\mathcal{T}_{\rm min}=\{E_{1}=\{1,3\}, E_{2}=\{7,8\}, E_{3}=\{2,4,6\}\}$$ and also the complex $\Gamma_{G}$ has one connected component which is a vertex, namely $\{E_{3}\}$, and one connected component which is an edge, namely $\{E_{1},E_{2}\}$.}
\end{ex1}

In \cite{K} the author studied the binomial arithmetical rank of $I_{\mathcal{A}_{G}}$ in two cases, namely when $G$ is bipartite or $I_{\mathcal{A}_{G}}$ is generated by quadratic binomials. Notice that in both cases the graph $G$ satisfies $(\sharp)$. Every bipartite graph satisfies this condition, since it has no odd cycles. Also, from Theorem 1.2 in \cite{OH1}, every graph $G$, such that $I_{\mathcal{A}_{G}}$ is generated by quadratic binomials, satisfies condition $(\sharp)$.

\begin{rem1} {\rm If $G$ is bipartite or $I_{\mathcal{A}_{G}}$ is generated by quadratic binomials, then the toric ideal $I_{\mathcal{A}_{G}}$ has the following property: $I_{\mathcal{A}_{G}}$ has no minimal binomials of the form $B_{w}=B_{w}^{(+)}-B_{w}^{(-)}$, where $B_{w}^{(+)}$, $B_{w}^{(-)}$ are squarefree monomials that are not indispensable of $I_{\mathcal{A}_{G}}$. To see this we distinguish the following cases: \begin{enumerate} \item the graph $G$ is bipartite. By Theorem 3.2 of \cite{Katz} the toric ideal $I_{\mathcal{A}_{G}}$ is minimally generated by all binomials of the form $B_{w}$, where $w$ is an even cycle with no chord. Now Theorem 3.2 of \cite{OH2} implies that every such binomial $B_{w}$ is indispensable. Thus every monomial arising in the unique minimal binomial generating set of $I_{\mathcal{A}_{G}}$ is indispensable.
\item the toric ideal $I_{\mathcal{A}_{G}}$ is generated by quadratic binomials. It is well known that the $\mathcal{A}_{G}$-degrees of the polynomials appearing in any minimal system of $\mathcal{A}_{G}$-homogeneous generators of $I_{\mathcal{A}_{G}}$ do not depend on the system of generators, see \cite[Section 8.3]{MS}. Using this fact and also that $I_{\mathcal{A}_{G}}$ has a quadratic set of binomial generators, we deduce that all minimal binomials of $I_{\mathcal{A}_{G}}$ are quadratic. Thus every monomial arising in a minimal system of binomial generators of $I_{\mathcal{A}_{G}}$ is quadratic and therefore, from Remark \ref{Extremal}, it is indispensable of $I_{\mathcal{A}_{G}}$.
    \end{enumerate}}
\end{rem1}

The next Theorem determines certain classes of toric ideals $I_{\mathcal{A}_{G}}$ for which the equality ${\rm bar}(I_{\mathcal{A}_{G}})=\mu(I_{\mathcal{A}_{G}})$ holds.

\begin{thm1} \label{BasicThm} Let $G$ be a graph such that $I_{\mathcal{A}_{G}}$ has no minimal binomials of the form $B_{w}=B_{w}^{(+)}-B_{w}^{(-)}$, where $B_{w}$ is a circuit and $B_{w}^{(+)}$, $B_{w}^{(-)}$ are squarefree monomials that are not indispensable of $I_{\mathcal{A}_{G}}$. Then ${\rm bar}(I_{\mathcal{A}_{G}})=\mu(I_{\mathcal{A}_{G}})$.
\end{thm1}

\noindent \noindent \textbf{Proof.} Since $G$ satisfies condition $(\sharp)$, the toric ideal $I_{\mathcal{A}_{G}}$ has a minimal binomial generating set $\mathcal{P}$ consisting only of circuits of the form $B_{w}$, where $B_{w}^{(+)}$ and $B_{w}^{(-)}$ are squarefree monomials. Notice that for each $B_{w} \in \mathcal{P}$ the walk $w$ is either an even cycle or two odd cycles intersecting in exactly one vertex. Given a binomial $B_{w} \in \mathcal{P}$, we have, from the assumption, that either exactly one of the monomials $B_{w}^{(+)}$, $B_{w}^{(-)}$ is indispensable or both of them are indispensable.\\ If $M$ is an indispensable monomial of $I_{\mathcal{A}_{G}}$, which is not quadratic, such that $\{E={\rm supp}(M)\}$ is a connected component of $\Gamma_{G}$, then there exists at least one binomial $B_{w}=M-{\bf x}^{\bf u} \in \mathcal{P}$ with the property that $M$ is a monomial of $B_{w}$. We will prove that $B_{w}$ is the unique binomial in $\mathcal{P}$ with the above property. Suppose that there exists another binomial $B_{\gamma} \in \mathcal{P}$ such that $B_{\gamma}=M-{\bf x}^{{\bf v}}$. Notice that the monomials ${\bf x}^{{\bf u}}$, ${\bf x}^{{\bf v}}$ are not indispensable, because $\{E\}$ is a connected component of $\Gamma_{G}$. Certainly $g={\bf x}^{{\bf u}}-{\bf x}^{{\bf v}}$ is a minimal binomial of $I_{\mathcal{A}_{G}}$ and therefore it is primitive. So ${\rm supp}({\bf x}^{{\bf u}}) \cap {\rm supp}({\bf x}^{{\bf v}})=\emptyset$ and also $g=B_{\zeta}$, for an even closed walk $\zeta$. Since the binomials $B_{w}$ and $B_{\gamma}$ belong to $\mathcal{P}$, the monomials ${\bf x}^{{\bf u}}$, ${\bf x}^{{\bf v}}$ are squarefree, so $\zeta$ is either an even cycle or two odd cycles intersecting in exactly one vertex. In fact the minimal binomial $B_{\zeta} \in I_{\mathcal{A}_{G}}$ is a circuit and it is a difference of two squarefree non-indispensable monomials, a contradiction to our assumption.\\
Let $q \geq 0$, $r \geq 0$ be the number of connected components of $\Gamma_{G}$ which are vertices and 2-simplices, correspondingly. Denote by $s=\mu(I_{\mathcal{A}_{G}})$ the minimal number of generators of $I_{\mathcal{A}_{G}}$, which is equal to the cardinality of the set $\mathcal{P}$, and also by $t \geq 0$ the number of indispensable binomials of $I_{\mathcal{A}_{G}}$. Proposition \ref{Edge} (2) asserts that $\Gamma_{G}$ has exactly $t$ connected components which are edges. Our aim is to prove that $r=\frac{s-q-t}{2}$. Let $B_{w_1}=M_{i}-M_{j} \in \mathcal{P}$ be a quadratic binomial that is not indispensable of $I_{\mathcal{A}_{G}}$, then the edge $\{E_{1}={\rm supp}(M_{i}),E_{2}={\rm supp}(M_{j})\}$ is not a connected component of $\Gamma_{G}$. Thus there exists an indispensable monomial $M_{k}$ of $I_{\mathcal{A}_{G}}$, such that $\{E_{1},E_{2},E_{3}={\rm supp}(M_{k})\}$ is a 2-simplex of $\Gamma_{G}$. Consider the quadratic binomials $B_{w_{2}}=M_{i}-M_{k} \in I_{\mathcal{A}_{G}}$, $B_{w_{3}}=M_{j}-M_{k} \in I_{\mathcal{A}_{G}}$ and notice that both of them are not indispensable of $I_{\mathcal{A}_{G}}$. Since all monomials $M_{i}$, $M_{j}$, $M_{k}$ are indispensable of $I_{\mathcal{A}_{G}}$ and $\{E_{1},E_{2},E_{3}\}$ is a 2-simplex of $\Gamma_{G}$, we deduce that there are exactly two binomials in $\mathcal{P}$ whose monomials are $M_{i}$, $M_{j}$ and $M_{k}$. Therefore $\Gamma_{G}$ has at least $\frac{s-q-t}{2}$ connected components which are 2-simplices, so $r \geq \frac{s-q-t}{2}$.\\
Let $E={\rm supp}(M_{i})$, $E'={\rm supp}(M_{j})$, $E''={\rm supp}(M_{k})$ be three elements of $\mathcal{T}_{\rm min}$ such that $\{E,E',E''\}$ is a 2-simplex of $\Gamma_{G}$, then the monomials $M_{i}$, $M_{j}$, $M_{k}$ are all of them at the same time quadratic and indispensable. Furthermore, there are quadratic binomials $B_{w_1}=M_{i}-M_{j} \in I_{\mathcal{A}_{G}}$, $B_{w_2}=M_{j}-M_{k} \in I_{\mathcal{A}_{G}}$, $B_{w_3}=M_{i}-M_{k} \in I_{\mathcal{A}_{G}}$. The minimal generating set $\mathcal{P}$ contains exactly two of the binomials $B_{w_1}$, $-B_{w_1}$, $B_{w_2}$, $-B_{w_2}$, $B_{w_3}$ and $-B_{w_3}$. Then $\mathcal{P}$ contains at least $2r$ binomials which are not indispensable. So $2r \leq s-q-t$. Consequently $r=\frac{s-q-t}{2}$.\\ For every connected component $\Gamma_{G}^{i}$ of $\Gamma_{G}$ which is a vertex we have that $\delta(\Gamma_{G}^{i})_{\{0,1\}}=1$, while for every connected component $\Gamma_{G}^{j}$ of $\Gamma_{G}$, which is an edge, we have that $\delta(\Gamma_{G}^{j})_{\{0,1\}}=1$. Also for every connected component $\Gamma_{G}^{i}$ of $\Gamma_{G}$ which is a 2-simplex, we have $\delta(\Gamma_{G}^{i})_{\{0,1\}}=2$. Consequently $$\delta(\Gamma_{G})_{\{0,1\}}=q+t+2 \frac{s-q-t}{2}=s;$$ i.e. $\delta(\Gamma_{G})_{\{0,1\}}=\mu(I_{\mathcal{A}_{G}})$, so, from Theorem \ref{Basic}, the inequality ${\rm bar}(I_{\mathcal{A}_{G}}) \geq \mu(I_{\mathcal{A}_{G}})$ holds and therefore ${\rm bar}(I_{\mathcal{A}_{G}})=\mu(I_{\mathcal{A}_{G}})$. \hfill $\square$\\

Theorem \ref{BasicThm} is no longer true if the toric ideal $I_{\mathcal{A}_{G}}$ has minimal binomials of the above form.
\begin{ex1}{\rm Consider the graph $G$ on the vertex set $\{1,\ldots,10\}$ with 14 edges, namely $e_{1}=\{1,2\}$, $e_{2}=\{2,3\}$, $e_{3}=\{3,4\}$, $e_{4}=\{4,5\}$, $e_{5}=\{5,6\}$, $e_{6}=\{6,7\}$, $e_{7}=\{7,8\}$, $e_{8}=\{8,9\}$, $e_{9}=\{9,10\}$, $e_{10}=\{1,10\}$, $e_{11}=\{1,5\}$, $e_{12}=\{2,6\}$, $e_{13}=\{1,7\}$ and $e_{14}=\{6,10\}$. The toric ideal $I_{\mathcal{A}_{G}}$ is minimally generated by the following nine binomials: $B_{w_1}=x_{1}x_{14}-x_{10}x_{12}$, $B_{w_2}=x_{5}x_{10}-x_{11}x_{14}$, $B_{w_3}=x_{6}x_{10}-x_{13}x_{14}$, $B_{w_4}=x_{1}x_{5}-x_{11}x_{12}$, $B_{w_5}= x_{5}x_{13}-x_{6}x_{11}$, $B_{w_6}=x_{1}x_{6}-x_{12}x_{13}$, $B_{w_7}=x_{2}x_{4}x_{6}x_{8}x_{14}-x_{3}x_{5}x_{7}x_{9}x_{12}$, $B_{w_8}=x_{1}x_{3}x_{7}x_{9}x_{11} - x_{2}x_{4}x_{8}x_{10}x_{13}$, $B_{w_9}=x_{1}x_{3}x_{5}x_{7}x_{9}-x_{2}x_{4}x_{6}x_{8}x_{10}$.\\ By Theorem \ref{Circuits} every binomial $B_{w_{i}}$ is a circuit, so $G$ satisfies condition $(\sharp)$. We have that the second power of $B_{w_9}$ belongs to the ideal generated by the binomials $B_{w_i}$, $1 \leq i \leq 8$, so ${\rm bar}(I_{\mathcal{A}_{G}}) \leq 8$. Using Theorem \ref{Basic} it is not hard to prove that ${\rm bar}(I_{\mathcal{A}_{G}}) \geq 8$, so in fact ${\rm bar}(I_{\mathcal{A}_{G}})=8$. Notice that the monomials $B_{w_9}^{(+)}=x_{1}x_{3}x_{5}x_{7}x_{9}$, $B_{w_9}^{(-)}=x_{2}x_{4}x_{6}x_{8}x_{10}$ are not indispensable. Actually the even cycle $w_{9}$ has two $\mathcal{F}_{4}$'s, namely $w_{3}=(e_{6},e_{13},e_{10},e_{14})$ and $w_{4}=(e_{1},e_{11},e_{5},e_{12})$.}
\end{ex1}

We prove now that the equality ${\rm bar}(I_{\mathcal{A}_{G}})=\mu(I_{\mathcal{A}_{G}})$ holds when $G$ is the wheel graph.
\begin{ex1} {\rm Consider the wheel graph $W_{n+1}$, $n \geq 3$, introduced in Example \ref{Wheel}. We will prove that ${\rm bar}(I_{\mathcal{A}_{W_{n+1}}})=\mu(I_{\mathcal{A}_{W_{n+1}}})$. If $n$ is even, then there exists, from Proposition 5.5. in \cite{OH3}, a bipartite graph $G$ such that $I_{\mathcal{A}_{W_{n+1}}}=I_{\mathcal{A}_{G}}$ and therefore we have, from Theorem 3.2 in \cite{K}, that ${\rm bar}(I_{\mathcal{A}_{W_{n+1}}})=\mu(I_{\mathcal{A}_{W_{n+1}}})$. If $n=3$, then $W_{4}$ is the complete graph on the vertex set $\{v_{1},\ldots,v_{4}\}$ and therefore $I_{\mathcal{A}_{W_4}}$ is complete intersection of height 2. Let us suppose that $n \geq 5$ is odd and assume that there is a minimal binomial $B_{w}=B_{w}^{(+)}-B_{w}^{(-)}$ of $I_{\mathcal{A}_{W_{n+1}}}$, where $B_{w}$ is a circuit, the monomials $B_{w}^{(+)}$, $B_{w}^{(-)}$ are squarefree and at least one of them is not indispensable of $I_{\mathcal{A}_{W_{n+1}}}$. Then the binomial $B_{w}$ is not indispensable and therefore, from Remark \ref{F4} (2), the walk $w$ has an $\mathcal{F}_{4}$, namely $\xi=(e_{1},e_{2},e_{3},e_{4})$ where $e_{2}$, $e_{4}$ are odd chords which cross effectively in $w$. Theorem 4.13 of \cite{RTT} implies that no odd chord of $w$ crosses the $\mathcal{F}_{4}$. By Definition \ref{Cross} the only possible case is $w=(\{v_{1},v_{2}\},\{v_{2},v_{3}\},\ldots,\{v_{n-1},v_{n}\},\{v_{n},v_{n+1}\},\{v_{1},v_{n+1}\})$, $e_{2}=\{v_{1},v_{n}\}$, $e_{4}=\{v_{2},v_{n+1}\}$ and $\xi=(e_{1}=\{v_{1},v_{2}\},e_{4},e_{3}=\{v_{n},v_{n+1}\},e_{2})$. Then $w$ can be written as $w=(w_{1},e_{3},w_{2},e_{1})$ where $w_{1}=(\{v_{2},v_{3}\},\{v_{3},v_{4}\},\ldots,\{v_{n-1},v_{n}\})$ and $w_{2}=\{v_{1},v_{n+1}\}$. But then the odd chord $\{v_{4},v_{n+1}\}$ of $w$ crosses the $\mathcal{F}_{4}$, a contradiction. By Theorem \ref{BasicThm} it holds that ${\rm bar}(I_{\mathcal{A}_{W_{n+1}}})=\mu(I_{\mathcal{A}_{W_{n+1}}})$.}
\end{ex1}

The {\em complement} of a graph $G$, denoted by $\overline{G}$, is the graph with the same vertices as $G$, and there is an edge between the vertices $v_{i}$ and $v_{j}$ if and only if there is no edge between $v_{i}$ and $v_{j}$ in $G$. A finite connected graph $G$ is called {\em weakly chordal} if every cycle of $G$ of length $4$ has a chord. In \cite{OH2} they study the toric ideal of a graph $G$ such that $\overline{G}$ is weakly chordal. We will prove that the equality ${\rm bar}(I_{\mathcal{A}_{G}})=\mu(I_{\mathcal{A}_{G}})$ holds for such graphs.

\begin{rem1} \label{Weakly} {\rm It follows easily that $\overline{G}$ is weakly chordal if and only if the following condition is satisfied: If $e$ and $e'$ are edges of $G$ with $e \cap e'=\emptyset$, then there is an edge $e''$ of $G$ with $e \cap e'' \neq \emptyset$ and also $e' \cap e'' \neq \emptyset$.}
\end{rem1}

\begin{prop1} Let $G$ be a graph such that $\overline{G}$ is weakly chordal, then ${\rm bar}(I_{\mathcal{A}_{G}})=\mu(I_{\mathcal{A}_{G}})$.
\end{prop1}
\noindent \textbf{Proof.} First we will prove that $G$ satisfies condition $(\sharp)$. Let $w$ be an even closed walk, which consists of two vertex disjoint odd cycles $\xi_{1}$ and $\xi_{2}$ joined by a path of length $\geq 1$ connecting one vertex $i$ of $\xi_{1}$ with one vertex $j$ of $\xi_{2}$. There are edges $e$, $e'$ of $G$ such that
\begin{enumerate} \item $e$ is an edge of $\xi_{1}$, which does not contain $i$ as a vertex.
\item $e'$ is an edge of $\xi_{2}$. \item $e \cap e'=\emptyset$.
\end{enumerate}
By Remark \ref{Weakly} there is an edge $e''$ of $G$ with $e \cap e'' \neq \emptyset$ and also $e' \cap e'' \neq \emptyset$. Thus $w$ can't be an induced subgraph of $G$.\\ Next we will prove that $I_{\mathcal{A}_{G}}$ has no minimal binomials of the form $B_{w}=B_{w}^{(+)}-B_{w}^{(-)}$, where $B_{w}$ is a circuit and $B_{w}^{(+)}$, $B_{w}^{(-)}$ are squarefree monomials that are not indispensable of $I_{\mathcal{A}_{G}}$. It follows from (Second Step) (a) of \cite[page 520]{OH1} that $I_{\mathcal{A}_{G}}$ has no minimal binomials of the form $B_{w}=B_{w}^{(+)}-B_{w}^{(-)}$, where $w$ is two odd cycles intersecting in exactly one vertex. In particular $I_{\mathcal{A}_{G}}$ is minimally generated by binomials of the form $B_{w}$, where $w$ is an even cycle. Let $B_{\xi}$ be a minimal binomial of $I_{\mathcal{A}_{G}}$, where $\xi=(e_{1}=\{1,2\},e_{2}=\{2,3\},\ldots,e_{2s}=\{2s,1\})$ is an even cycle of $G$ of length $2s \geq 6$. We will show that $\xi$ has no two odd chords which cross effectively. Let $e=\{1,2i+1\}$, $e'=\{2j,2k\}$ be two odd chords which cross effectively in $\xi$, i.e. $1<2j<2i+1<2k \leq 2s$. Since $B_{\xi}$ is a minimal binomial, we have, from Theorem 4.13 of \cite{RTT}, that the chords $e$, $e'$ can't cross strongly effectively and therefore they form an $\mathcal{F}_{4}$, denoted by $\gamma$, of $\xi$. The above theorem implies that $\xi$ has no even chords and also that there is no odd chord of $\xi$ which crosses the $\mathcal{F}_{4}$. Consider the edges $e_{2}=\{2,3\}$ and $e_{2s-1}=\{2s-1,2s\}$ and notice that they share no common vertex. Thus there exists an edge $e_{t}$ of $G$ with $e_{t} \cap e_{2} \neq \emptyset$ and $e_{t} \cap e_{2s-1} \neq \emptyset$. Certainly $e_{t}$ is a chord of $\xi$. Using again the fact that $B_{\xi}$ is a minimal generator of $I_{\mathcal{A}_{G}}$, we have, from Theorem 4.13 in \cite{RTT}, that $e_{t}$ is odd chord. Thus either $e_{t}=\{2,2s\}$ or $e_{t}=\{3,2s-1\}$. Let us assume that $k \neq s$. Then necessarily $\{1,2j\}$ is an edge of $\gamma$, so it is an edge of $\xi$ and therefore $j=1$. Also $\{2i+1,2k\}$ is an edge of $\xi$, so $2k=2i+2$ and therefore $k=i+1$. Thus $\gamma=(e_{1},e',e_{2i+1}=\{2i+1,2i+2\},e)$. The even cycle $\xi$ can be written as $\xi=(\xi_{1},e_{1},\xi_{2},e_{2i+1})$, where $\xi_{1}=(e_{2k}=\{2k,2k+1\},e_{2k+1},\ldots,e_{2s})$ and $\xi_{2}=(e_{2},e_{3},\ldots,e_{2i})$. Now the odd chord $e_{t}$ crosses the $\mathcal{F}_{4}$, a contradiction. Assume now that $k=s$ and let $j \neq 1$. Then $\{1,2k\}$ is an edge of $\xi$ and also $\{2j,2i+1\}$ is an edge of $\xi$, so $2i+1=2j+1$ and therefore $i=j$. Thus $\gamma=(e_{2s},e',e_{2j}=\{2i,2i+1\},e)$. The even cycle $\xi$ can be written as $\xi=(\xi_{1},e_{2s},\xi_{2},e_{2i})$, where $\xi_{1}=(e_{2i+1},e_{2i+2},\ldots,e_{2s-1})$ and $\xi_{2}=(e_{1},e_{2},\ldots,e_{2i-1})$. The odd chord $e_{t}$ of $\xi$ crosses the $\mathcal{F}_{4}$, a contradiction. We work analogously for the case $j=1$ and arrive at a contradiction. Thus $\xi$ has no two odd chords which cross effectively, so we have, from Theorem 3.2 of \cite{OH2}, that $B_{\xi}$ is indispensable of $I_{\mathcal{A}_{G}}$ and therefore the monomials $B_{\xi}^{(+)}$, $B_{\xi}^{(-)}$ are indispensable. Now Theorem \ref{BasicThm} implies that ${\rm bar}(I_{\mathcal{A}_{G}})=\mu(I_{\mathcal{A}_{G}})$. \hfill $\square$\\

\section{A class of determinantal lattice ideals}

Let $S=K[x_{1},\ldots,x_{m},y_{1},\ldots,y_{m}]$ be the polynomial ring in $2m$ variables with coefficients in a field $K$. Consider the ideal $I_{2}(D) \subset S$ generated by the 2-minors of the $2 \times m$ matrix of indeterminants $$D=\begin{pmatrix}
x_{1}^{d_{1}} & x_{2}^{d_{2}} & \ldots & x_{m}^{d_{m}}\\
y_{1}^{d_{1}} & y_{2}^{d_{2}} & \ldots & y_{m}^{d_{m}}
\end{pmatrix},$$where every $d_{i}$, $1 \leq i \leq m$, is a positive integer. When $d_{1}=d_{2}=\cdots=d_{m}=1$ the quotient $S/I_{2}(D)$ is the coordinate ring of the Segre embedding $\mathbb{P}_{K}^{1} \times \mathbb{P}_{K}^{m}$. For  $1 \leq i<j \leq m$ we let $f_{ij}:=x_{i}^{d_i}y_{j}^{d_j}-x_{j}^{d_j}y_{i}^{d_i}$.

\begin{thm1} \label{Grobner} The reduced Gr{\"o}bner basis with respect to any term order $\prec$ in $S$ for the ideal $I_{2}(D)$ is given by $\mathcal{G}=\{f_{ij}|1 \leq i<j \leq m\}$.
\end{thm1}
\noindent \textbf{Proof.} Consider two binomials $f_{ij} \in \mathcal{G}$ and $f_{kl} \in \mathcal{G}$. We will prove that $S(f_{ij},f_{kl}) \stackrel{\mathcal{G}}{\longrightarrow} 0$. Let us first examine the case that ${\rm in}_{\prec}(f_{ij})=x_{i}^{d_{i}}y_{j}^{d_{j}}$ and ${\rm in}_{\prec}(f_{kl})=x_{k}^{d_{k}}y_{l}^{d_{l}}$. If $i \neq k$ and $j \neq l$, then $S(f_{ij},f_{kl}) \stackrel{\mathcal{G}}{\longrightarrow} 0$ since the initial monomials are relatively prime. Suppose that $i=k$, so $j \neq l$. Without loss of generality we can assume that $j<l$. We have that $S(f_{ij},f_{kl})=x_{l}^{d_{l}}y_{j}^{d_{j}}y_{k}^{d_{k}}-x_{j}^{d_{j}}y_{l}^{d_{l}}y_{k}^{d_{k}} \stackrel{f_{jl}}{\longrightarrow} 0$. Let $j=l$, then $i \neq k$. Without loss of generality we can assume that $i>k$. We have that $S(f_{ij},f_{kl})=x_{i}^{d_{i}}y_{k}^{d_{k}}x_{j}^{d_{j}}-x_{k}^{d_{k}}y_{i}^{d_{i}}x_{j}^{d_{j}} \stackrel{f_{ki}}{\longrightarrow} 0$. Using similar arguments we take that $S(f_{ij},f_{kl}) \stackrel{\mathcal{G}}{\longrightarrow} 0$ in the remaining cases, namely the cases \begin{enumerate} \item ${\rm in}_{\prec}(f_{ij})=x_{i}^{d_{i}}y_{j}^{d_{j}}$ and ${\rm in}_{\prec}(f_{kl})=x_{l}^{d_{l}}y_{k}^{d_{k}}$. \item ${\rm in}_{\prec}(f_{ij})=x_{j}^{d_{j}}y_{i}^{d_{i}}$ and ${\rm in}_{\prec}(f_{kl})=x_{k}^{d_{k}}y_{l}^{d_{l}}$.
\item ${\rm in}_{\prec}(f_{ij})=x_{j}^{d_{j}}y_{i}^{d_{i}}$ and ${\rm in}_{\prec}(f_{kl})=x_{l}^{d_{l}}y_{k}^{d_{k}}$.
\end{enumerate}
Consequently $\mathcal{G}$ is a Gr\"obner basis of $I_{2}(D)$, with respect to any term order $\prec$. Clearly it is also a reduced Gr\"obner basis of $I_{2}(D)$. \hfill $\square$\\

\begin{rem1} {\rm It is clear that $\mathcal{G}$ is a minimal generating set of the ideal $I_{2}(D)$.}
\end{rem1}

\begin{prop1} \label{Lattice} The ideal $I_{2}(D)$ is a lattice ideal of height $m-1$.
\end{prop1}
\noindent \textbf{Proof.} By Theorem \ref{Grobner} the set $\mathcal{G}$ is the reduced Gr\"obner basis of $I_{2}(D)$ with respect to the graded reverse lexicographic term order induced by any ordering of the variables $x_{i}$ and $y_{i}$, $1 \leq i \leq m$. Lemma 12.1 of \cite{St} applies and guarantees that $I_{2}(D):x_{i}^{\infty}=I_{2}(D)$ and $I_{2}(D):y_{i}^{\infty}=I_{2}(D)$, for every $1 \leq i \leq m$. Thus $$(I_{2}(D):(x_{1} \cdots x_{m})^{\infty})=((((I_{2}(D):x_{1}^{\infty}):x_{2}^{\infty}) \cdots):x_{m}^{\infty})=I_{2}(D)$$ and similarly $(I_{2}(D):(y_{1} \cdots y_{m})^{\infty})=I_{2}(D)$. Therefore $$(I_{2}(D):(x_{1} \cdots x_{m}y_{1} \cdots y_{m})^{\infty})=I_{2}(D),$$ so we deduce, from Corollary 2.5 of \cite{E-S}, that $I_{2}(D)$ is a lattice ideal. As a consequence $I_{2}(D)$ is of the form $I_{L}$, for a lattice $L \subset \mathbb{Z}^{2m}$. For $1 \leq i<j \leq m$ we let ${\bf u}_{ij} \in \mathbb{Z}^{m}$ the vector with coordinates $$({\bf u}_{ij})_{k}=\begin{cases} d_{i}, & \textrm{if} \  k=i \\ -d_{j}, & \textrm{if} \  k=j\\
0, &  \textrm{otherwise}.
\end{cases}$$Since $\mathcal{G}$ is a set of generators for $I_{2}(D)$, we have, from Lemma 2.5 in \cite{LV}, that the set of all vectors ${\bf v}_{ij}:=({\bf u}_{ij},-{\bf u}_{ij}) \in \mathbb{Z}^{2m}$, $1 \leq i<j \leq m$, generates the lattice $L$. Furthermore, for every $1<i<j \leq m$ we have that ${\bf v}_{ij}={\bf v}_{1j}-{\bf v}_{1i}$. Thus $L$ is generated by all the vectors ${\bf v}_{1j}=({\bf u}_{1j},-{\bf u}_{1j})$, $2 \leq j \leq m$. Since all the above vectors are $\mathbb{Z}$-linearly independent, we have that ${\rm rank}(L)=m-1$. Thus ${\rm ht}(I_{2}(D))=m-1$. \hfill $\square$\\

\begin{rem1} \label{RemarkGrobner} {\rm Consider the monomial ideals $\mathcal{M}_{1}=(x_{1}^{d_{1}}y_{j}^{d_j}|2 \leq j \leq m)$ and $\mathcal{M}_{2}=(x_{i}^{d_i}y_{j}^{d_j}| 2 \leq i<j \leq m)$. Then the initial ideal of $I_{2}(D)$ with respect to the lexicographic term order $\prec$ induced by $x_{1} \succ x_{2} \succ \cdots \succ x_{m} \succ y_{1} \succ  \cdots \succ y_{m}$ is equal to the sum $\mathcal{M}_{1}+\mathcal{M}_{2}$.}
\end{rem1}

\begin{not1} {\rm For the rest of this section we will keep the notation introduced in Theorem \ref{Grobner} and Proposition \ref{Lattice}.}
\end{not1}

\begin{prop1} \label{IndispensableLattice} The ideal $I_{2}(D)$ is generated by its indispensable.
\end{prop1}
\noindent \textbf{Proof.} By Remark \ref{RemIndi} the set $\{x_{i}^{d_{i}}y_{j}^{d_{j}},x_{j}^{d_{j}}y_{i}^{d_{i}}| 1 \leq i<j \leq m\}$ generates the monomial ideal $\mathcal{M}_{L}$. Actually it is a minimal generating set for the ideal $\mathcal{M}_{L}$, so for every $1 \leq i<j \leq m$ the monomials $x_{i}^{d_{i}}y_{j}^{d_{j}}$, $x_{j}^{d_{j}}y_{i}^{d_{i}}$ are indispensable of $I_{2}(D)$. Let $\mathcal{P}=\{B_{1},\ldots,B_{s}\}$ be a minimal binomial generating set of $I_{2}(D)$. The monomial $x_{i}^{d_{i}}y_{j}^{d_{j}}:={\bf x}^{\bf u}$ is indispensable of $I_{2}(D)$, so there exists $1 \leq k \leq s$ such that $B_{k}={\bf x}^{\bf u}-{\bf x}^{\bf v}$. Our aim is to prove that $B_{k}=f_{ij}$. Notice that none of the variables $x_{i}$ and $y_{j}$ divides ${\bf x}^{\bf v}$. If at least one of the variables $x_{i}$, $y_{j}$ divides ${\bf x}^{\bf v}$, then ${\bf x}^{\bf w}={\rm gcd}({\bf x}^{\bf u},{\bf x}^{\bf v}) \neq 1$ and therefore the binomial ${\bf x}^{{\bf u}-{\bf w}}-{\bf x}^{{\bf v}-{\bf w}} \in I_{2}(D)$. Thus ${\bf x}^{{\bf u}-{\bf w}} \in \mathcal{M}_{L}$ and properly divides ${\bf x}^{\bf u}$, a contradiction to the fact that ${\bf x}^{\bf u}$ is a minimal generator of $\mathcal{M}_{L}$. Assume that ${\bf x}^{\bf v} \neq x_{j}^{d_{j}}y_{i}^{d_{i}}$ and consider now the non-zero binomial $g=x_{j}^{d_{j}}y_{i}^{d_{i}}-{\bf x}^{\bf v}$. Since the monomial $x_{j}^{d_{j}}y_{i}^{d_{i}}$ is indispensable of $I_{2}(D)$, using similar arguments as before we take that none of the variables $x_{j}$ and $y_{i}$ divides ${\bf x}^{\bf v}$. By Theorem \ref{Grobner} we have that $B_{k} \stackrel{\mathcal{G}}{\longrightarrow} 0$, because $B_{k} \in I_{2}(D)$, and this can happen only if ${\bf x}^{\bf v}=x_{j}^{d_{j}}y_{i}^{d_{i}}$. Thus $B_{k}=f_{ij}$. Analogously it can be proved that $f_{ij}$ is the only binomial in $\mathcal{P}$ which contains $x_{j}^{d_{j}}y_{i}^{d_{i}}$ as a monomial. \hfill $\square$\\

\begin{rem1} \label{Clear} {\rm Clearly the support of any indispensable monomial of $I_{2}(D)$ belongs to $\mathcal{T}_{min}$. Thus $\Gamma_{L}$ has exactly $m(m-1)$ vertices.}
\end{rem1}

The next theorem asserts that the equality ${\rm bar}(I_{2}(D))=\mu(I_{2}(D))$ holds.

\begin{thm1} For the lattice ideal $I_{2}(D)$ we have ${\rm bar}(I_{2}(D))=\frac{m(m-1)}{2}$.
\end{thm1}
\noindent \textbf{Proof.} By Remark \ref{Clear} the simplicial complex $\Gamma_{L}$ has $m(m-1)$ vertices and therefore, from Proposition \ref{Integer}, ${\rm bar}(I_{2}(D)) \geq \frac{m(m-1)}{2}$. Since $\mu(I_{2}(D))$ equals $\frac{m(m-1)}{2}$, we have that also ${\rm bar}(I_{2}(D)) \leq \frac{m(m-1)}{2}$, so ${\rm bar}(I_{2}(D))=\frac{m(m-1)}{2}$. \hfill $\square$\\

Consider the vector configuration $\mathcal{B}=\{b_{1},\ldots,b_{m}\} \subset \mathbb{N}$, for  $$b_{i}=\prod_{\substack{j=1\\ j \neq i}}^{m}d_{j}, \ 1 \leq i \leq m.$$ Let $\mathcal{A}$ be the set of columns of the Lawrence lifting of $\mathcal{B}$, namely the $(m+1) \times 2m$-matrix $$\Lambda(\mathcal{B})=\begin{pmatrix}
\mathcal{B} & {\bf 0}_{1 \times m}\\
{\bf I}_{m} & {\bf I}_{m}
\end{pmatrix}$$ where ${\bf I}_{m}$ is the $m \times m$ identity matrix and ${\bf 0}_{1 \times m}$ is the $1 \times m$ zero matrix. The toric ideal $I_{\mathcal{A}} \subset S$ is the kernel of the $K$-algebra homomorphism $\phi:S \rightarrow K[t_{1},\ldots,t_{m+1}]$ given by $\phi(x_{i})=t_{1}^{b_{i}}t_{i+1}$ and $\phi(y_{i})=t_{i+1}$, for $1 \leq i \leq m$. From the definitions $I_{2}(D)$ is contained in the toric ideal $I_{\mathcal{A}}$ which has height $m-1$. By \cite[Corollary 2.2]{E-S} $I_{Sat(L)}$ is the only minimal prime of $I_{L}$ which is a binomial ideal. Thus $Sat(L)=ker_{\mathbb{Z}}(\mathcal{A})$. For every $1 \leq i \leq m$ we let $d_{i}^{\star}=\frac{d_{i}}{{\rm gcd}(d_{i},d_{j})}$.

\begin{prop1} The ideal $I_{2}(D)$ is prime if and only if it holds that ${\rm gcd}(d_{i},d_{j})=1$, for every $1 \leq i<j \leq m$.
\end{prop1}
\noindent \textbf{Proof.} Let us assume that $I_{2}(D)$ is prime, i.e. $I_{2}(D)=I_{\mathcal{A}}$, and also that there exist $1 \leq i<j \leq m$ such that ${\rm gcd}(d_{i},d_{j}) \neq 1$. Since the vector ${\bf v}_{ij}$ belongs to $L$, we have that $\frac{1}{{\rm gcd}(d_{i},d_{j})} {\bf v}_{ij}$ belongs to $Sat(L)=ker_{\mathbb{Z}}(\mathcal{A})$ and therefore the binomial $x_{i}^{d_{i}^{\star}}y_{j}^{d_{j}^{\star}}-x_{j}^{d_{j}^{\star}}y_{i}^{d_{i}^{\star}}$ belongs to $I_{2}(D)$. Thus $x_{i}^{d_{i}^{\star}}y_{j}^{d_{j}^{\star}} \in \mathcal{M}_{L}$ and properly divides $x_{i}^{d_{i}}y_{j}^{d_{j}}$, a contradiction to the fact that $x_{i}^{d_{i}}y_{j}^{d_{j}}$ is a minimal generator of $\mathcal{M}_{L}$. Conversely assume that ${\rm gcd}(d_{i},d_{j})=1$, for every $1
\leq i<j \leq m$. Notice that ${\rm gcd}(b_{1},\ldots,b_{m})=1$. We have, from Proposition 10.1.8 of \cite{Vil1}, that the lattice $ker_{\mathbb{Z}}(\mathcal{B})$ is generated by all vectors ${\bf u}_{ij}$, $1 \leq i<j \leq m$, and therefore $ker_{\mathbb{Z}}(\mathcal{A})$ is generated by all vectors ${\bf v}_{ij}$. Thus $L=ker_{\mathbb{Z}}(\mathcal{A})$ and therefore $I_{2}(D)$ is prime. \hfill $\square$\\

The following theorem provides a lower bound for the binomial arithmetical rank of $I_{\mathcal{A}}$.

\begin{thm1} \label{Lawrence} For the toric ideal $I_{\mathcal{A}}$ we have ${\rm bar}(I_{\mathcal{A}}) \geq \frac{m(m-1)}{2}$.
\end{thm1}
\noindent \textbf{Proof.} Given a monomial ${\bf x}^{\bf v}{\bf y}^{\bf w} \in S$, we let ${\rm supp}({\bf x}^{\bf v}{\bf y}^{\bf w})={\rm supp}({\bf z})$, for ${\bf z}=({\bf v},{\bf w}) \in \mathbb{N}^{2m}$. First we show that $I_{\mathcal{A}}$ contains no binomials of the form $x_{i}^{u_{i}}-{\bf x}^{\bf v}{\bf y}^{\bf w}$ or $y_{j}^{u_{j}}-{\bf x}^{\bf v}{\bf y}^{\bf w}$. Let, say, that $I_{\mathcal{A}}$ has a binomial $B=x_{i}^{u_{i}}-{\bf x}^{\bf v}{\bf y}^{\bf w}$. Since $Sat(L)=ker_{\mathbb{Z}}(\mathcal{A})$, there exists a positive integer $d$ such that $x_{i}^{du_{i}}-{\bf x}^{d{\bf v}}{\bf y}^{d{\bf w}} \in I_{2}(D)$. Thus $x_{i}^{du_{i}} \in \mathcal{M}_{L}$ and therefore it should be divided by an indispensable monomial of $I_{2}(D)$, a contradiction. Similarly we can prove that $I_{\mathcal{A}}$ has no binomials of the form $y_{j}^{u_{j}}-{\bf x}^{\bf v}{\bf y}^{\bf w}$.\\ For every $1 \leq i<j \leq m$ we have that both monomials $x_{i}^{d_{i}}y_{j}^{d_j}$, $x_{j}^{d_{j}}y_{i}^{d_i}$ are in the monomial ideal $\mathcal{M}_{ker_{\mathbb{Z}}(\mathcal{A})}$, so the ideal $\mathcal{M}_{ker_{\mathbb{Z}}(\mathcal{A})}$ has two minimal generators $M_{ij}$, $M_{ji}$ such that ${\rm supp}(M_{ij})={\rm supp}(x_{i}^{d_{i}}y_{j}^{d_j})$ and ${\rm supp}(M_{ji})={\rm supp}(x_{j}^{d_{j}}y_{i}^{d_i})$. By Proposition 1.5 of \cite{KO}, the monomials $M_{ij}$, $M_{ji}$ are indispensable of $I_{\mathcal{A}}$. Also their support is minimal with respect to inclusion, i.e. there exists no monomial $N \in \mathcal{M}_{ker_{\mathbb{Z}}(\mathcal{A})}$ with ${\rm supp}(N) \subsetneqq M_{ij}$ or ${\rm supp}(N) \subsetneqq M_{ji}$. Thus $\Gamma_{ker_{\mathbb{Z}}(\mathcal{A})}$ has at least $m(m-1)$ vertices and therefore we have, from Proposition \ref{Integer}, that ${\rm bar}(I_{\mathcal{A}}) \geq \frac{m(m-1)}{2}$. \hfill $\square$\\

We give now an example of a toric ideal $I_{\mathcal{A}}$ such that ${\rm bar}(I_{\mathcal{A}})=\frac{m(m-1)}{2}$.

\begin{ex1} {\rm Let $d_{1}=2$, $d_{2}=4$, $d_{3}=5$ and $d_{4}=7$. Then $b_{1}=140$, $b_{2}=70$, $b_{3}=56$ and $b_{4}=40$. Thus $\mathcal{A}$ is the set of columns of the matrix $$\begin{pmatrix}
140 & 70 & 56 & 40 & 0 & 0 & 0 & 0
            \\
            1 & 0 & 0 & 0 & 1 & 0 & 0 & 0
            \\
            0 & 1 & 0 & 0 & 0 & 1 & 0 & 0
            \\
            0 & 0 & 1 & 0 & 0 & 0 & 1 & 0
            \\
            0 & 0 & 0 & 1 & 0 & 0 & 0 & 1
            \end{pmatrix}.$$ The toric ideal $I_{\mathcal{A}}$ is minimally generated by the following 8 binomials:\\ $B_{1}=x_{1}y_{2}^{2}-x_{2}^{2}y_{1}$, $B_{2}=x_{1}^{2}y_{3}^{5}-x_{3}^{5}y_{1}^2$, $B_{3}=x_{1}^2y_{4}^7-x_{4}^7y_{1}^2$, $B_{4}=x_{2}^4y_{3}^5-x_{3}^{5}y_{2}^{4}$, $B_{5}=x_{2}^4y_{4}^7-x_{4}^7y_{2}^4$, $B_{6}=x_{3}^5y_{4}^7-x_{4}^7y_{3}^5$, $B_{7}=x_{1}x_{2}^2y_{3}^5-x_{3}^5y_{1}y_{2}^2 $, $B_{8}=x_{4}^7y_{1}y_{2}^2-x_{1}x_{2}^2y_{4}^7$.\\ Here $m=4$, so Theorem \ref{Lawrence} implies that ${\rm bar}(I_{\mathcal{A}}) \geq 6$. Furthermore, we have that $rad(I_{\mathcal{A}})=rad(B_{1},\ldots,B_{6})$, since the second power of $B_{7}$, as well as the second power of $B_{8}$, belongs to the ideal generated by the binomials $B_{1},\ldots,B_{6}$. Thus ${\rm bar}(I_{\mathcal{A}})=6$.}
\end{ex1}

As it was proved in Proposition \ref{Lattice} the set $\{{\bf v}_{1j}| 2 \leq j \leq m\}$ is a $\mathbb{Z}$-basis for the lattice $L$. Let $J_{L} \subset S$ be the ideal generated by all binomials $f_{1j}=x_{1}^{d_{1}}y_{j}^{d_{j}}-x_{j}^{d_{j}}y_{1}^{d_{1}}$ where $2 \leq j \leq m$. The ideal $J_{L}$ is commonly known as the {\em lattice basis ideal} of $L$. We will compute the minimal primary decomposition of $rad(J_{L})$, when $I_{2}(D)$ is prime.

\begin{lem1} \label{Grobner1} The set $\mathcal{R}=\{f_{12},f_{13},\ldots,f_{1m}\} \cup \{g_{ij}:=y_{1}^{d_{1}}f_{ij}|2 \leq i<j \leq m\}$ is a Gr{\"o}bner basis of the lattice basis ideal $J_{L}$ with respect to the lexicographic term order $\prec$ induced by $x_{1} \succ x_{2} \succ \cdots \succ x_{m} \succ y_{1} \succ  \cdots \succ y_{m}$.
\end{lem1}
\noindent \textbf{Proof.} First we prove that $\mathcal{R} \subset J_{L}$. Since $\{f_{1k}| 2 \leq k \leq m\} \subset J_{L}$, it is enough to show that $g_{ij} \in J_{L}$. For every $2 \leq i<j \leq m$ we have that $$g_{ij}=y_{1}^{d_{1}}(x_{i}^{d_{i}}y_{j}^{d_{j}}-x_{j}^{d_{j}}y_{i}^{d_{i}})=y_{i}^{d_{i}}f_{1j}-y_{j}^{d_{j}}f_{1i} \in J_{L}.$$ Let $f_{1k}=x_{1}^{d_{1}}y_{k}^{d_{k}}-x_{k}^{d_{k}}y_{1}^{d_{1}}$ and $f_{1l}=x_{1}^{d_{1}}y_{l}^{d_{l}}-x_{l}^{d_{l}}y_{1}^{d_{1}}$, where $2 \leq k<l \leq m$. It holds that $$S(f_{1k},f_{1l})=y_{1}^{d_{1}}x_{l}^{d_l}y_{k}^{d_k}-y_{1}^{d_{1}}x_{k}^{d_k}y_{l}^{d_l} \stackrel{g_{kl}}{\longrightarrow} 0.$$ We will prove that $S(f_{1k},g_{ij}) \stackrel{\mathcal{R}}{\longrightarrow} 0$. If $j \neq k$, then the initial monomials ${\rm in}_{\prec}(f_{1k})=x_{1}^{d_{1}}y_{k}^{d_{k}}$, ${\rm in}_{\prec}(g_{ij})=x_{i}^{d_{i}}y_{1}^{d_{1}}y_{j}^{d_{j}}$ are relatively prime and therefore $S(f_{1k},g_{ij}) \stackrel{\mathcal{R}}{\longrightarrow} 0$. If $j=k$, then $$S(f_{1k},g_{ij})=x_{1}^{d_{1}}x_{k}^{d_{k}}y_{1}^{d_{1}}y_{i}^{d_{i}}-x_{i}^{d_{i}}y_{1}^{d_{1}}x_{k}^{d_{k}}y_{1}^{d_{1}} \stackrel{f_{1i}}{\longrightarrow} 0.  \ \ \ \square$$

\begin{prop1} \label{BasicPro} For the lattice basis ideal $J_{L}$ we have $J_{L}=I_{2}(D) \cap (x_{1}^{d_{1}},y_{1}^{d_{1}})$.
\end{prop1}
\noindent \textbf{Proof.} Clearly $J_{L} \subset I_{2}(D) \cap \mathcal{Q}$, where $\mathcal{Q}=(x_{1}^{d_{1}},y_{1}^{d_{1}})$. It remains to prove that $I_{2}(D) \cap \mathcal{Q} \subset J_{L}$. Consider the lexicographic term order $\prec$ in $S$ induced by $x_{1} \succ x_{2} \succ \cdots \succ x_{m} \succ y_{1} \succ \cdots \succ y_{m}$. Since $${\rm in}_{\prec}(J_{L}) \subset {\rm in}_{\prec} (I_{2}(D) \cap \mathcal{Q}) \subset {\rm in}_{\prec}(I_{2}(D)) \cap {\rm in}_{\prec}(\mathcal{Q}),$$ it is enough to prove that ${\rm in}_{\prec}(I_{2}(D)) \cap {\rm in}_{\prec}(\mathcal{Q}) \subset {\rm in}_{\prec}(J_{L})$. Notice that ${\rm in}_{\prec}(\mathcal{Q})=(x_{1}^{d_{1}},y_{1}^{d_{1}})=\mathcal{Q}$ and also $\{x_{1}^{d_{1}},y_{1}^{d_{1}}\}$ is a minimal generating set of $\mathcal{Q}$. Consider the monomial ideals $\mathcal{M}_{1}=(x_{1}^{d_{1}}y_{j}^{d_j}|2 \leq j \leq m)$ and $\mathcal{M}_{2}=(x_{i}^{d_i}y_{j}^{d_j}| 2 \leq i<j \leq m)$. Lemma \ref{Grobner1} asserts that the set $\mathcal{R}$ is a Gr{\"o}bner basis of the lattice basis ideal $J_{L}$ with respect to $\prec$, so ${\rm in}_{\prec}(J_{L})=\mathcal{M}_{1}+y_{1}^{d_{1}}\mathcal{M}_{2}$. We have, from Remark \ref{RemarkGrobner}, that ${\rm in}_{\prec}(I_{2}(D))$ is equal to the sum $\mathcal{M}_{1}+\mathcal{M}_{2}$. Actually $\{x_{1}^{d_{1}}y_{j}^{d_j}|2 \leq j \leq m\} \cup \{x_{i}^{d_i}y_{j}^{d_j}| 2 \leq i<j \leq m\}$ is a minimal generating set of ${\rm in}_{\prec}(I_{2}(D))$. For every $2 \leq j \leq m$ it holds that ${\rm lcm}(x_{1}^{d_{1}}y_{j}^{d_{j}},x_{1}^{d_{1}})=x_{1}^{d_{1}}y_{j}^{d_{j}} \in \mathcal{M}_{1} \subset {\rm in}_{\prec}(J_{L})$ and ${\rm lcm}(x_{1}^{d_{1}}y_{j}^{d_{j}},y_{1}^{d_{1}})=y_{1}^{d_{1}}x_{1}^{d_{1}}y_{j}^{d_{j}} \in \mathcal{M}_{1}$. Furthermore, for every $2 \leq i<j \leq m$ we have that ${\rm lcm}(x_{i}^{d_{i}}y_{j}^{d_{j}},x_{1}^{d_{1}})=x_{i}^{d_{i}}x_{1}^{d_{1}}y_{j}^{d_{j}} \in \mathcal{M}_{1} $ and ${\rm lcm}(x_{i}^{d_{i}}y_{j}^{d_{j}},y_{1}^{d_{1}})=y_{1}^{d_{1}}x_{i}^{d_{i}}y_{j}^{d_{j}} \in y_{1}^{d_{1}}\mathcal{M}_{2} \subset {\rm in}_{\prec}(J_{L})$. Proposition 1.2.1 of \cite{HH} implies that $${\rm in}_{\prec}(I_{2}(D)) \cap \mathcal{Q} \subset {\rm in}_{\prec}(J_{L}). \ \ \ \square$$

\begin{thm1} Suppose that for every $1 \leq i<j \leq m$ it holds that ${\rm gcd}(d_{i},d_{j})=1$. Then the minimal primary decomposition of the radical of the ideal $J_{L}$ is $$rad(J_{L})=I_{2}(D) \cap (x_{1},y_{1}).$$
\end{thm1}
\noindent \textbf{Proof.} Using Proposition \ref{BasicPro} we have that $$rad(J_{L})=rad(I_{2}(D)) \cap rad(x_{1}^{d_{1}},y_{1}^{d_{1}})=I_{2}(D) \cap (x_{1},y_{1}). \ \ \ \square $$


\begin{thebibliography}{00}

\bibitem{CKT} H. Charalambous, A. Katsabekis and A. Thoma, Minimal systems of binomial generators and the indispensable complex of a toric ideal, {\em Proc. Amer. Math. Soc.} {\bf 135} (2007) 3443-3451.
\bibitem{E-S} D. Eisenbud and B. Sturmfels, Binomial ideals, {\em Duke Math. J.} \textbf{84} (1996) 1-45.
\bibitem{HH} J. Herzog and T. Hibi, {\em Monomial ideals}, Graduate Texts in Mathematics {\bf 260}, Springer, 2010.
\bibitem{K} A. Katsabekis, Arithmetical rank of toric ideals associated to graphs, {\em Proc. Amer. Math. Soc.} {\bf 138} (2010) 3111-3123.
    \bibitem{KMT} A.Katsabekis, M. Morales and A. Thoma, Stanley-Reisner rings and the radicals of lattice ideals, {\em J. Pure Appl. Algebra} \textbf{204} (2006) 584-601.
         \bibitem{KO} A. Katsabekis and I. Ojeda, An indispensable classification of monomial curves in $\mathbb{A}^{4}(K)$, arXiv:1103.4702.
    \bibitem{KT} A. Katsabekis and A. Thoma, Matchings in simplicial complexes, circuits and toric varieties, {\em J. Comb. Theory Ser. A} {\bf 114} (2007) 300-310.
               \bibitem{Katz} M. Katzmann, Bipartite graphs whose edge algebras are complete intersections, {\em J. Algebra} {\bf 220} (1999) 519-530.
\bibitem{LV} H. L{\'o}peza and R. Villarreal, Complete intersections in binomial and lattice ideals, arXiv: 1205.0772.
    \bibitem{MS} E. Miller and B. Sturmfels, {\em Combinatorial Commutative Algebra}, Vol 227 of Graduate Texts in Mathematics, Springer, New York, 2005.
    \bibitem{OH1} H. Ohsugi and T. Hibi, Toric ideals generated by quadratic binomials, {\em J. Algebra} {\bf 218} (1999) 509-527.
        \bibitem{OH2} H. Ohsugi and T. Hibi, Indispensable binomials of toric ideals, {\em J. Algebra Appl.} {\bf 4} (2005) 421-434.
           \bibitem{OH3} H. Ohsugi and T. Hibi, Centrally symmetric configurations of integer matrices, {\em Nagoya Mathematical Journal} (to appear).
\bibitem{OH} H. Ohsugi and T. Hibi, Toric ideals and their circuits, {\em Journal of Commutative Algebra} (to appear).
    \bibitem{RTT} E. Reyes, C. Tatakis and A. Thoma, Minimal generators of toric ideals of graphs, {\em Adv. in Appl. Math.} {\bf 48} (2012) 64-78.
        \bibitem{St} B. Sturmfels, {\em Gr{\"o}bner Bases and Convex Polytopes}, University Lecture Series, No. 8 American Mathematical Society Providence, R.I. 1995.
            \bibitem{Vil} R. Villarreal, Rees algebras of edge ideals, {\em Comm. Algebra} {\bf 23} (1995) 3513-3524.
            \bibitem{Vil1} R. Villarreal, {\em Monomial Algebras}, Pure Appl. Math., Marcel Dekker, New York, 2001.


\end{thebibliography}
\end{document}